\documentclass[11pt,a4paper]{article}

\usepackage{inputenc}
\usepackage{amsmath}
\usepackage{bm}
\usepackage{bbold}
\usepackage{amsthm}
\usepackage{enumerate}

\usepackage{algorithm}

\usepackage[top=1.0in, left=1.0in, right=1.0in, bottom=1.0in]{geometry}

\usepackage{tikz}\tikzset{x=1cm,y=1cm,z=1cm}

\usepackage{pgfplots}\pgfplotsset{compat=1.16}

\usepackage[hyphens]{url}

\usepackage{hyperref}
\usepackage{breakurl}

\title{On solution of tropical discrete best approximation problems\thanks{Soft Computing 28(20), 12097-12112 (2024).}}

\author{N. Krivulin\thanks{Faculty of Mathematics and Mechanics, St.~Petersburg State University, 28 Universitetsky Ave., St.~Petersburg, 198504, Russia; 
nkk@math.spbu.ru.}
}

\date{}

\newtheorem{theorem}{theorem}
\newtheorem{lemma}[theorem]{lemma}

\theoremstyle{definition}

\begin{document}

\maketitle

\begin{abstract}
We consider a discrete best approximation problem formulated in the framework of tropical algebra, which deals with the theory and applications of algebraic systems with idempotent operations. Given a set of samples of input and output of an unknown function, the problem is to construct a generalized tropical Puiseux polynomial that best approximates the function in the sense of a tropical distance function. The construction of an approximate polynomial involves the evaluation of both unknown coefficient and exponent of each monomial in the polynomial. To solve the approximation problem, we first reduce the problem to an equation in unknown vector of coefficients, which is given by a matrix with entries parameterized by unknown exponents. We derive a best approximate solution of the equation, which yields both vector of coefficients and approximation error parameterized by the exponents. Optimal values of exponents are found by minimization of the approximation error, which is transformed into minimization of a function of exponents over all partitions of a finite set. We solve this minimization problem in terms of max-plus algebra (where addition is defined as maximum and multiplication as arithmetic addition) by using a computational procedure based on the agglomerative clustering technique. This solution is extended to the minimization problem of finding optimal exponents in the polynomial in terms of max-algebra (where addition is defined as maximum). The results obtained are applied to develop new solutions for conventional problems of discrete best Chebyshev approximation of real functions by piecewise linear functions and piecewise Puiseux polynomials. We discuss computational complexity of the proposed solution and estimate upper bounds on the computational time. We demonstrate examples of approximation problems solved in terms of max-plus and max-algebra, and give graphical illustrations.
\\

\textbf{Key-Words:} tropical semifield, tropical Puiseux polynomial, best approximate solution, discrete best approximation, Chebyshev approximation.
\\

\textbf{MSC (2020):} 15A80, 90C24, 41A50, 41A65, 65D15
\end{abstract}

\section{Introduction}

Discrete best approximation in which an unknown function is approximated from sample data appears in a variety of theoretical and applied contexts. The problems to reconstruct a function from sample observations as accurate as possible are frequently encountered in many research areas where the uncertainty in the data, model or acquisition technique is the key issue to be addressed, including interval, fuzzy and neural computing \cite{Dubois1980Fuzzy,Samarasinghe2006Neural,Celikyilmaz2009Modeling,Schunn2010Uncertainly}.

Given a set of points and corresponding values of the function, the goal is to find a continuous function of a particular form, which minimizes approximation error in the sense of a certain metric. In many cases, the approximating functions are assumed to be piecewise polynomials (including piecewise linear functions), whereas the error is measured by the Chebyshev metric, which leads to Chebyshev best approximation problems \cite{Mhaskar2000Fundamentals,Steffens2006History} that date back to Laplace's classical work \cite{Laplace1832Mecanique} (Book 3, Chapter V, \S39).

To address the discrete best approximation problems with both Chebyshev and other metrics, a variety of methods were proposed in the literature (see, e.g., overviews in \cite{Conn1988Computational,Szusz2010Linear}). Available techniques are based on linear programming \cite{Osborne1967Best,Watson1970Algorithm}, dynamic programming \cite{Gluss1962Further,Camponogara2015Models}, and other optimization approaches \cite{Stone1961Approximation,Cameron1966Piecewise,Tomek1974Two,Imai1986Optimal,Szusz2010Linear}. These techniques can provide optimal or near-optimal solutions with a moderate polynomial computational complexity.

Discrete best approximation problems can be formulated and solved in the framework of tropical algebra \cite{Golan2003Semirings,Gondran2008Graphs,Heidergott2006Maxplus,Butkovic2010Maxlinear,Maclagan2015Introduction}, which deals with the theory of semirings and semifields with idempotent addition and finds applications in a range of areas from algebraic geometry to operations research. As an example of idempotent semifields, one can consider max-plus algebra where addition is defined as maximum and multiplication as arithmetic addition. Another example is max-algebra where addition is defined as maximum and multiplication as usual. The application of tropical algebra allows many problems, which are nonlinear in the ordinary sense, to be transformed into linear problems in tropical algebra, and hence facilitates the formal analysis and simplifies the derivation of solution. 

The implementation of tropical algebra to approximation problems mainly focuses on the best approximation of a vector in a tropical linear space \cite{Butkovic2010Maxlinear,Akian2011Best,Saad2021Zerosum}. The discrete best approximation problems of functions, which appear in recent works on neural networks and machine learning (see, e.g., \cite{Zhang2018Tropical,Maragos2021Tropical}), constitute another direction of tropical approximation that needs further research. As an attempt to address this need, a general problem of discrete best approximation of a function is introduced and examined in tropical algebra settings in \cite{Krivulin2023Algebraic}, where a function defined on a tropical semifield is approximated by tropical Puiseux polynomials and rational functions with respect to a generalized metric.

The tropical Puiseux polynomials \cite{Markwig2010Field,Grigoriev2018Tropical}, which take the same form as in the conventional algebra with the difference that the addition and multiplication are defined in terms of tropical algebra, have applications in a variety of fields, including algebraic geometry \cite{Itenberg2007Tropical}, image processing \cite{Li1992Morphological}, cryptography \cite{Grigoriev2014Tropical,Grigoriev2019Tropical} and games \cite{Esparza2008Approximative}. Each Puiseux polynomial written in terms of max-plus algebra is represented as a piecewise linear function in the standard notation, whereas the tropical metric coincides with the standard Chebyshev metric. When defined in the framework of max-algebra, the Puiseux polynomials correspond to ordinary piecewise Puiseux polynomial functions, which may be of interest in spline approximation. These correspondences offer a potential to solve conventional best Chebyshev approximation problems through their transformation into tropical approximation problems, which makes the development of new efficient solutions of tropical discrete best approximation problems highly relevant.

In this paper, we consider a discrete best approximation problem formulated in the tropical algebra setting. Given a set of samples of corresponding input and output values of an unknown function, the problem is to construct a generalized tropical Puiseux polynomial (a polynomial with real exponents) that best approximates the function from the sample data. 

We further develop the approach described in \cite{Krivulin2023Algebraic}, where the exponents (powers) are assumed to be fixed in advance, and the solution concentrates on evaluating the coefficients of the monomials in the approximating polynomial. Since the situation when the exponents are also unknown is much more common, we now solve a more general problem where the construction of an approximating polynomial involves the evaluation of both unknown coefficient and exponent of each monomial.

To solve the approximation problem, we first reduce it to a vector equation in the unknown vector of coefficients, which is given by a matrix with entries parameterized by the unknown exponents. We derive a best approximate solution of the equation by using the approach developed in \cite{Krivulin2009Solution,Krivulin2012Solution,Krivulin2013Solution-linear}, which results in both vector of coefficients and approximation error parameterized by the exponents. Optimal values of exponents are then found by minimization of the approximation error, which is transformed into minimization of an objective function over all partitions of a finite set. For each partition, this function is evaluated by solving a set of minimization problems each corresponding to a part of the partition. The maximum value of minimums of these problems is taken as the value of the objective function, whereas the solutions of the problems determine the exponents.

We solve the above minimization problem in the framework of max-plus algebra, where the evaluation of the objective function is reduced to the solution of a set of polynomial optimization problems in which the unknown exponents become indeterminates. To find an optimal (near-optimal) partition, we propose a computational procedure based on the agglomerative clustering technique. The procedure starts with the partition into single element subsets and then forms better partitions by subsequent merging subsets by pairs. For each subset in the partitions, the procedure solves its related polynomial optimization problem by applying the solution proposed in \cite{Krivulin2021Algebraic}. 

The exponents obtained by the procedure are used to evaluate both the corresponding vector of coefficients and approximation error, which completes the derivation of the approximating polynomial. Furthermore, this solution technique is extended to the minimization problem of finding optimal exponents in the approximating polynomial in terms of max-algebra. The results obtained are applied to develop new solutions for conventional problems of discrete best Chebyshev approximation of real functions by piecewise linear functions and piecewise Puiseux polynomials. We discuss computational complexity of the proposed solution and estimate upper bounds on the computational time required. We demonstrate examples of approximation problems solved in terms of max-plus and max-algebra, and give graphical illustrations.

The rest of the paper is organized as follows. In Section~\ref{S-PDNR}, we present main definitions, notation and results of tropical algebra, which we use in what follows. Section~\ref{S-PAF} introduces the polynomial approximation problem under study and describes the transformation of the problem to a form suitable for further analysis and solution. A solution procedure for the approximation in terms of max-plus and max-algebra is proposed in Section~\ref{S-AMPMAP}. We offer numerical examples in Section~\ref{S-INE} and give concluding remarks in Section~\ref{S-C}.

\section{Preliminary definitions, notation and results}
\label{S-PDNR}

In this section, we outline the key definitions and notation, and describe some preliminary results of tropical (idempotent) algebra, which provide an analytical framework for the solutions developed in subsequent sections. For further details on the theory and methods of tropical mathematics, one can consult a range of textbooks and monographs published in recent decades, including \cite{Golan2003Semirings,Gondran2008Graphs,Heidergott2006Maxplus,Butkovic2010Maxlinear,Maclagan2015Introduction}.

\subsection{Tropical semifield}

Let $\mathbb{X}$ be a set that is closed under addition $\oplus$ and multiplication $\otimes$, and equipped with zero $\mathbb{0}$ and identity $\mathbb{1}$. We assume that $(\mathbb{X},\mathbb{0},\oplus)$ is a commutative idempotent monoid, $(\mathbb{X}\setminus\{\mathbb{0}\},\mathbb{1},\otimes)$ is an abelian group, and multiplication distributes over addition. The algebraic structure $(\mathbb{X},\mathbb{0},\mathbb{1},\oplus,\otimes)$ is called a tropical (or idempotent) semifield.

In the semifield, addition is idempotent, which implies that $x\oplus x=x$ for each $x\in\mathbb{X}$. Multiplication is invertible to ensure that each $x\ne\mathbb{0}$ has its inverse $x^{-1}$ such that $x\otimes x^{-1}=\mathbb{1}$. Hereafter, we drop the multiplication sign $\otimes$ for the sake of brevity.

Throughout the paper (unless explicitly stated otherwise), we assume the same algebraic context defined in the setting of a general tropical semifield, where the abbreviation $xy$ is interpreted as the tropical product $x\otimes y$. To avoid possible misinterpretation of the abbreviation in the cases where a different meaning is implied such as arithmetic multiplication, we give a remark to indicate directly the standard arithmetic context or use explicit formula $x\times y$.

The power notation with integer exponents serves to represent iterated products defined as $x^{p}=xx^{p-1}$, $x^{-p}=(x^{-1})^p$, $x^{0}=\mathbb{1}$ and $\mathbb{0}^{p}=\mathbb{0}$ for any nonzero $x\in\mathbb{X}$ and positive integer $p$. It is assumed that the equation $x^{p}=a$ has a unique solution $x$ for any $a\in\mathbb{X}$ and integer $p>0$, which gives rise to the powers with rational exponents. Moreover, we consider that the rational powers can be further extended to allow real exponents. In what follows, the powers are thought of as defined in terms of tropical algebra unless otherwise indicated.

Idempotent addition generates a partial order on $\mathbb{X}$ by the rule: $x\leq y$ if and only if $x\oplus y=y$. We assume the condition $x\oplus y\in\{x,y\}$ holds for all $x,y$, which makes the semifield selective and hence turns the above partial order into a total order. With respect to this order, both addition and multiplication are monotone, which means that the inequality $x\leq y$ results in $x\oplus z\leq y\oplus z$ and $xz\leq yz$ for any $z\in\mathbb{X}$. Furthermore, addition has an extremal property (the majority law), which says that $x\leq x\oplus y$ and $y\leq x\oplus y$. The inequality $x\oplus y\leq z$ is equivalent to the system of inequalities $x\leq z$ and $y\leq z$. Finally, it follows from the inequality $x\leq y$ for some nonzero $x,y\in\mathbb{X}$, that $x^{r}\geq y^{r}$ if $r<0$ and $x^{r}\leq y^{r}$ if $r\geq0$.

According to the extremal property of addition, we can define a maximum operation (with respect to the order induced by idempotent addition) as $\max(x,y)=x\oplus y$ for any $x,y\in\mathbb{X}$. As it is not difficult to see, a dual minimum operation can then be defined as $\min(x,y)=(x^{-1}\oplus y^{-1})^{-1}$ if $x,y\ne\mathbb{0}$, and $\min(x,y)=\mathbb{0}$ otherwise.

As examples of the idempotent semifield under consideration, one can consider the real semifields
\begin{equation*}
\mathbb{R}_{\max,+}
=
(\mathbb{R}\cup\{-\infty\},-\infty,0,\max,+),
\qquad
\mathbb{R}_{\max}
=
(\mathbb{R}_{+},0,1,\max,\times),
\end{equation*}
where $\mathbb{R}$ is the set of reals, and $\mathbb{R}_{+}=\{x\in\mathbb{R}|\ x\geq0\}$.

In the semifield $\mathbb{R}_{\max,+}$, which is commonly referred to as max-plus algebra, the addition $\oplus$ is defined as $\max$ and multiplication $\otimes$ as $+$. The zero $\mathbb{0}$ is given by $-\infty$ and identity $\mathbb{1}$ by $0$. For each $x\in\mathbb{R}$, the inverse $x^{-1}$ corresponds to the opposite number $-x$ in regular notation. The power $x^{y}$ coincides with the conventional product $x\times y$. The order induced by the idempotent addition agrees with the natural linear order on $\mathbb{R}$.

The semifield $\mathbb{R}_{\max}$, also known as max-algebra, has $\oplus=\max$, $\otimes=\times$, $\mathbb{0}=0$ and $\mathbb{1}=1$. The inverse and power notations are interpreted as usual. The order associated with addition is the standard linear order on $\mathbb{R}_{+}$.

Both idempotent semifields are isomorphic to each other by the following isomorphisms: $\exp:\mathbb{R}_{\max,+}\rightarrow\mathbb{R}_{\max}$ and $\log:\mathbb{R}_{\max}\rightarrow\mathbb{R}_{\max,+}$.

\subsection{Distributivity of addition and minimum}

We now present distributive identities for addition (maximum) and minimum to be used in the subsequent analysis. We start with an identity that represents the distributivity of $\oplus$ over $\min$. First note that for any $x,y,u,v\in\mathbb{X}$, the following equality obviously holds:
\begin{equation*}
\min(x,y)\oplus\min(u,v)
=
\min(x\oplus u,x\oplus v,y\oplus u,y\oplus v).
\end{equation*}

To rewrite this identity for an arbitrary number of terms, consider scalars $x_{ij}\in\mathbb{X}$ for all $i=1,\ldots,M$ and $j=1,\ldots,N$, where $M$ and $N$ are positive integers. By extending the above identity, we immediately obtain
\begin{equation*}
\bigoplus_{i=1}^{M}
\min_{1\leq j\leq N}
x_{ij}
=
\min_{1\leq j_{1},\ldots,j_{M}\leq N}
(x_{1j_{1}}\oplus\cdots\oplus x_{Mj_{M}}).
\end{equation*}

We observe that the summands $x_{ij}$ under the minimum operator on the right are in the ascending order of the first index. To rearrange the summands in the order of the second index, we consider the set of all ordered partitions $(I_{1},\ldots,I_{N})$ of the set $\{1,\ldots,M\}$ into $N$ parts (including empty parts).

After rearrangement by grouping summands with common second indices, the right-hand side of the last identity is rewritten as
\begin{equation*}
\min_{1\leq j_{1},\ldots,j_{M}\leq N}
(x_{1j_{1}}\oplus\cdots\oplus x_{Mj_{M}})
=
\min_{(I_{1},\ldots,I_{N})}
\left(
\bigoplus_{i\in I_{1}}
x_{i1}
\oplus\cdots\oplus
\bigoplus_{i\in I_{N}}
x_{iN}
\right),
\end{equation*}
where the right minimum is taken over all ordered set partitions $(I_{1},\ldots,I_{N})$. (Here and hereafter, we interpret the empty sums as $\mathbb{0}$.)

By combining the results, we arrive at the distributive identity
\begin{equation}
\bigoplus_{i=1}^{M}
\min_{1\leq j\leq N}
x_{ij}
=
\min_{(I_{1},\ldots,I_{N})}
\bigoplus_{j=1}^{N}
\bigoplus_{i\in I_{j}}
x_{ij}.
\label{E-plusmin}
\end{equation}

Next we consider a function $F(x_{1},\ldots,x_{N})=f_{1}(x_{1})\oplus\cdots\oplus f_{N}(x_{N})$, where $f_{j}:\mathbb{X}\rightarrow\mathbb{X}$ are given functions for all $j=1,\ldots,N$. The function $F(x_{1},\ldots,x_{N})$ takes the form of max-separable functions that were extensively studied in \cite{Zimmermann1984Maxseparable,Zimmermann2003Disjunctive,Tharwat2010One}. Due to the separability, the evaluation of the minimum of $F(x_{1},\ldots,x_{N})$ over all $x_{1},\dots,x_{N}$ reduces to minimization of each function $f_{j}(x_{j})$, which yields the following distributive identity: 
\begin{equation}
\min_{x_{1},\ldots,x_{N}}
\bigoplus_{j=1}^{N}
f_{j}(x_{j})
=
\bigoplus_{j=1}^{N}
\min_{x_{j}}
f_{j}(x_{j}).
\label{E-minplus}
\end{equation}

\subsection{Puiseux polynomials}

In this subsection, we concentrate on polynomial functions defined in the tropical algebra setting. We consider a generalized Puiseux polynomial in one variable $x$, given in the form
\begin{equation}
P(x)
=
\bigoplus_{j=1}^{N}
a_{j}
x^{p_{j}},
\qquad
x\ne\mathbb{0},
\label{E-ajxpj}
\end{equation}
where $a_{j}\ne\mathbb{0}$ are coefficients and $p_{j}\in\mathbb{R}$ are exponents for all $j=1,\ldots,N$.

If the polynomial is thought of in the context of max-plus algebra (the semifield $\mathbb{R}_{\max,+}$), it is given in terms of the usual operations as
\begin{equation*}
\max_{1\leq j\leq N}
(p_{j}\times x+a_{j}),
\end{equation*} 
which is a piecewise-linear convex function on $\mathbb{R}$.

For the polynomial defined in the context of max-algebra (the semifield $\mathbb{R}_{\max}$), the conventional form becomes the piecewise polynomial
\begin{equation*}
\max_{1\leq j\leq N}
a_{j}\times
x^{p_{j}},
\end{equation*} 
which may be either a convex or nonconvex function on $\mathbb{R}_{+}$.

Consider a max-algebra polynomial $P(x)$ given in the conventional form and take a logarithm to a base greater than $1$. It follows from the monotonicity of logarithm that
\begin{equation*}
\log P(x)
=
\log
\left(\max_{1\leq j\leq N}
a_{j}\times
x^{p_{j}}\right)
=
\max_{1\leq j\leq N}
(\log a_{j}+p_{j}\times\log x).
\end{equation*} 

This means that taking the logarithm maps $P(x)$ into a max-plus algebra polynomial
\begin{equation*}
P^{\prime}(x)
=
\max_{1\leq j\leq N}
(p_{j}\times x^{\prime}+a_{j}^{\prime}),
\end{equation*} 
where $x^{\prime}=\log x$ and $a_{j}^{\prime}=\log a_{j}$ for all $j=1,\ldots,N$.

Suppose that given polynomial \eqref{E-ajxpj}, one needs to solve a polynomial optimization problem to find $x\in\mathbb{X}$ that attains the minimum
\begin{equation}
\begin{aligned}
&&
\min_{x}
&&&
\bigoplus_{j=1}^{N}
a_{j}x^{p_{j}}.
\end{aligned}
\label{P-minxajxpj}
\end{equation}

The next result provides a solution as a special case of the result in \cite{Krivulin2021Algebraic}, obtained for the corresponding constrained optimization problem.
\begin{lemma}
\label{L-minxajxpj}
The minimum value in problem \eqref{P-minxajxpj} is equal to
\begin{equation}
\mu
=
\bigoplus_{\substack{1\leq j,k\leq N\\p_{j}<0,\ p_{k}>0}}
a_{j}^{-\frac{p_{k}}{p_{j}-p_{k}}}a_{k}^{\frac{p_{j}}{p_{j}-p_{k}}}
\oplus
\bigoplus_{\substack{1\leq j\leq N\\p_{j}=0}}
a_{j},
\label{E-mu}
\end{equation}
whereas all solutions are given by the condition
\begin{equation*}
\bigoplus_{\substack{1\leq j\leq N\\p_{j}<0}}
\mu^{1/p_{j}}
a_{j}^{-1/p_{j}}
\leq
x
\leq
\Biggl(
\bigoplus_{\substack{1\leq j\leq N\\p_{j}>0}}
\mu^{-1/p_{j}}
a_{j}^{1/p_{j}}
\Biggr)^{-1}.
\end{equation*}
\end{lemma}

We observe that the right-hand side of the double inequality can be rewritten as minimum to bring the inequality into the form
\begin{equation}
\bigoplus_{\substack{1\leq j\leq N\\p_{j}<0}}
\mu^{1/p_{j}}
a_{j}^{-1/p_{j}}
\leq
x
\leq
\min_{\substack{1\leq j\leq N\\p_{j}>0}}
\mu^{1/p_{j}}
a_{j}^{-1/p_{j}}.
\label{I-x}
\end{equation}

As it is not difficult to see, the computational complexity of this solution grows not faster than $O(N^{2})$.

\subsection{Matrix and vector algebra}

An idempotent semifield of scalars is extended to matrices over the semifield in the usual way. We denote by $\mathbb{X}^{M\times N}$ the set of matrices of $M$ rows and $N$ columns with entries from $\mathbb{X}$. A matrix that has all entries equal to $\mathbb{0}$ is the zero matrix. If a matrix has no zero rows (columns), it is called row-regular (column-regular). A matrix that is both row- and column-regular is called regular.

Matrix addition, matrix multiplication and scalar multiplication of matrix follow the standard rules where the scalar addition and multiplication are replaced by $\oplus$ and $\otimes$. For any conforming matrices $\bm{A}=(a_{ij})$, $\bm{B}=(b_{ij})$, $\bm{C}=(c_{ij})$, and scalar $x$, these matrix operations are performed according to the entrywise formulas
\begin{equation*}
(\bm{A}\oplus\bm{B})_{ij}
=
a_{ij}\oplus b_{ij},
\qquad
(\bm{A}\bm{C})_{ij}
=
\bigoplus_{k}a_{ik}c_{kj},
\qquad
(x\bm{A})_{ij}
=
xa_{ij}.
\end{equation*}

Monotonic properties of the scalar operations $\oplus$ and $\otimes$ extend to the matrix operations, where the inequalities are interpreted entrywise.

A matrix of one row (column) is a row (column) vector. The set of all column vectors with $N$ components is denoted $\mathbb{X}^{N}$. All vectors are considered column vectors unless otherwise specified. A vector that has all components equal to $\mathbb{0}$ is the zero vector denoted $\bm{0}$. If a vector has no zero components, it is called regular.

For any nonzero column vector $\bm{x}=(x_{i})$, its (multiplicative) conjugate is the row vector $\bm{x}^{-}=(x_{i}^{-})$ where $x_{i}^{-}=x_{i}^{-1}$ if $x_{i}\ne\mathbb{0}$, and $x_{i}^{-}=\mathbb{0}$ otherwise.

For any vector $\bm{x}=(x_{i})$ in $\mathbb{X}^{N}$, the support is defined as the set of indices of nonzero components in $\bm{x}$, given by $\mathop\mathrm{supp}(\bm{x})=\{i|\ x_{i}\ne\mathbb{0},\ 1\leq i\leq N\}$. Let $\bm{x}=(x_{i})$ and $\bm{y}=(y_{i})$ be nonzero vectors such that $\mathop\mathrm{supp}(\bm{x})=\mathop\mathrm{supp}(\bm{y})$. Define the distance between the vectors by the following distance function: 
\begin{equation}
d(\bm{x},\bm{y})
=
\bigoplus_{i\in\mathop\mathrm{supp}(\bm{x})}\left(x_{i}y_{i}^{-1}\oplus x_{i}^{-1}y_{i}\right)
=
\bm{y}^{-}\bm{x}
\oplus
\bm{x}^{-}\bm{y}.
\label{E-dxy}
\end{equation}

Under the condition $\mathop\mathrm{supp}(\bm{x})\ne\mathop\mathrm{supp}(\bm{y})$, we put $d(\bm{x},\bm{y})=\infty$, where $\infty$ denotes a value that is greater than any $x\in\mathbb{X}$ (an undefined value). Finally, if $\bm{x}=\bm{y}=\bm{0}$, then we assume $d(\bm{x},\bm{y})=\mathbb{1}$.

We observe that in the framework of the semifield $\mathbb{R}_{\max,+}$ (max-plus algebra) with $\mathbb{1}=0$, the distance function $d$ coincides for all $\bm{x},\bm{y}\in\mathbb{R}^{N}$ with the Chebyshev metric, which in terms of the regular algebra is given by
\begin{equation*}
d_{\infty}(\bm{x},\bm{y})
=
\max_{1\leq i\leq N}|x_{i}-y_{i}|
=
\max_{1\leq i\leq N}\max(x_{i}-y_{i},y_{i}-x_{i}).
\end{equation*}

In the case of the semifield $\mathbb{R}_{\max}$ (max-algebra), the function $d$ at \eqref{E-dxy} can be considered as a generalized metric that takes values in the interval $[1,\infty)\subset\mathbb{R}_{+}$. Moreover, taking into account the isomorphism between $\mathbb{R}_{\max}$ and $\mathbb{R}_{\max,+}$, this function transforms into a metric $d^{\prime}(\bm{x},\bm{y})=\log d(\bm{x},\bm{y})$.

\subsection{Solution of vector equation}

Given a matrix $\bm{A}\in\mathbb{X}^{M\times N}$ and vector $\bm{b}\in\mathbb{X}^{M}$, consider the problem to find regular vectors $\bm{x}\in\mathbb{X}^{N}$ that solve the equation
\begin{equation}
\bm{A}\bm{x}
=
\bm{b}.
\label{E-Axeqb}
\end{equation}

A vector $\bm{x}_{\ast}$ is a best approximate solution of equation \eqref{E-Axeqb} in the sense of the metric $d$,  if for all vectors $\bm{x}$ the following inequality holds:
\begin{equation*}
d(\bm{A}\bm{x}_{\ast},\bm{b})
\leq
d(\bm{A}\bm{x},\bm{b}).
\end{equation*}

The best approximate solution is achieved by finding
\begin{equation*}
\bm{x}_{\ast}
=
\arg\min_{\bm{x}}
d(\bm{A}\bm{x},\bm{b}).
\end{equation*}

The value of the distance $d(\bm{A}\bm{x}_{\ast},\bm{b})$ represents the approximation error associated with the best approximate solution.

A best approximate solution to equation \eqref{E-Axeqb} under general assumptions can be obtained following the results in \cite{Krivulin2009Solution,Krivulin2012Solution,Krivulin2013Solution-linear} in the next form.
\begin{theorem}
\label{T-Axeqb}
Let $\bm{A}$ be a regular matrix and $\bm{b}$ a regular vector. Define the scalar $\Delta=(\bm{A}(\bm{b}^{-}\bm{A})^{-})^{-}\bm{b}$. Then, the following statements hold:
\begin{enumerate}
\item
The best approximation error for equation \eqref{E-Axeqb} is equal to
\begin{equation*}
d(\bm{A}\bm{x}_{\ast},\bm{b})
=
\sqrt{\Delta};
\end{equation*}
\item
The best approximate solution of the equation is given by
\begin{equation*}
\bm{x}_{\ast}
=
\sqrt{\Delta}(\bm{b}^{-}\bm{A})^{-}.
\end{equation*}
\item
If $\Delta=\mathbb{1}$, then the equation has exact solutions, which include the maximal solution given by $\bm{x}_{\ast}=(\bm{b}^{-}\bm{A})^{-}$.
\end{enumerate}
\end{theorem}

\section{Polynomial approximation of functions}
\label{S-PAF}

We are now in a position to describe and analyze the tropical polynomial approximation problem under study. The problem aims to construct a best approximating generalized Puiseux polynomial by finding both coefficients and exponents in the monomials of the polynomial. To handle the problem, we propose an approach that involves two main steps. The first step of the solution consists in the transformation of the problem into a best approximation of a vector in a tropical vector space. This vector approximation problem is solved analytically, which results in both the coefficients and approximation error parameterized by the exponents.  

As the second step, we need to find exponents that minimize the approximation error. We use distributive identities for addition and minimum to represent the problem of minimizing the error as a minimization of an objective function over partitions of a finite set. Each part of a partition defines a function of one exponent that has to be minimized. The maximum value over all functions defined for parts is taken as the value of the objective function corresponding to the entire partition. The set of exponents that minimize the functions of parts is considered as optimal exponents that lead to the minimum approximation error.

Consider a problem of approximating a function $f:\mathbb{X}\rightarrow\mathbb{X}$ from nonzero samples $x_{i},y_{i}\in\mathbb{X}$, where $x_{i}$ and $y_{i}$ are an input and its corresponding output of the function, given for each $i=1,\ldots,M$. As an approximating function, we use a generalized Puiseux polynomial given a number $N$ of monomials.

For each $j=1,\ldots,N$, let $p_{j}\in\mathbb{R}$ denote an unknown exponent (power), and $\theta_{j}\in\mathbb{X}$ an unknown coefficient (parameter). The approximating polynomial is defined as
\begin{equation*}
P(x)
=
\bigoplus_{j=1}^{N}\theta_{j}x^{p_{j}}
=
\theta_{1}x^{p_{1}}
\oplus\cdots\oplus
\theta_{N}x^{p_{N}}.
\end{equation*}

The approximation problem is formulated to determine both unknown exponents and parameters from the input-output data by minimizing the deviation between the left and right sides of the equations
\begin{equation}
\theta_{1}x_{i}^{p_{1}}
\oplus\cdots\oplus
\theta_{N}x_{i}^{p_{N}}
=
y_{i},
\qquad
i=1,\ldots,M.
\label{E-theta1xip1thetaNxipN-yi}
\end{equation}

\subsection{Vector form of approximation problem}

To represent the system of equations at \eqref{E-theta1xip1thetaNxipN-yi} in vector form, we follow the approach in \cite{Krivulin2023Algebraic} to introduce the following vector-matrix notation:
\begin{gather*}
\bm{p}
=
\left(
\begin{array}{c}
p_{1}
\\
\vdots
\\
p_{N}
\end{array}
\right),
\qquad
\bm{\theta}
=
\left(
\begin{array}{c}
\theta_{1}
\\
\vdots
\\
\theta_{N}
\end{array}
\right),
\qquad
\bm{y}
=
\left(
\begin{array}{c}
y_{1},
\\
\vdots
\\
y_{M}
\end{array}
\right),
\\
\bm{X}(\bm{p})
=
\left(
\begin{array}{ccc}
x_{1}^{p_{1}} & \ldots & x_{1}^{p_{N}}
\\
\vdots & & \vdots
\\
x_{M}^{p_{1}} & \ldots & x_{M}^{p_{N}}
\end{array}
\right).
\end{gather*}

With this notation, the system combines into the vector equation
\begin{equation}
\bm{X}(\bm{p})\bm{\theta}
=
\bm{y}.
\label{E-Xpthetaeqy}
\end{equation}

The polynomial approximation problem now reduces to finding a best approximate solution in the sense of a minimum distance between the vectors $\bm{X}(\bm{p})\bm{\theta}$ and $\bm{y}$. With the distance function $d$, the problem takes the form
\begin{equation}
\begin{aligned}
&&
\min_{(\bm{\theta},\bm{p})}
&&&
d(\bm{X}(\bm{p})\bm{\theta},\bm{y}).
\end{aligned}
\label{P-mindXpthetay}
\end{equation}

To obtain a best approximation polynomial, one needs to find vectors
\begin{equation*}
\bm{\theta}_{\ast}
=
(\theta_{1}^{\ast},\ldots,\theta_{N}^{\ast})^{T},
\qquad
\bm{p}_{\ast}
=
(p_{1}^{\ast},\ldots,p_{N}^{\ast})^{T},
\end{equation*}
which solve the minimization problem at \eqref{P-mindXpthetay} as the minimizers
\begin{equation*}
(\bm{\theta}_{\ast},\bm{p}_{\ast})
=
\arg\min_{(\bm{\theta},\bm{p})}
d(\bm{X}(p)\bm{\theta},\bm{y}).
\end{equation*}

The best approximating polynomial is then defined to be
\begin{equation*}
P_{\ast}(x)
=
\theta_{1}^{\ast}x^{p_{1}^{\ast}}
\oplus\cdots\oplus
\theta_{N}^{\ast}x^{p_{N}^{\ast}}.
\end{equation*}

\subsection{Solution approach to approximation problem}

The purpose of this subsection is to outline a solution approach to be implemented to solve problem \eqref{P-mindXpthetay}. First, we split the minimum over pairs of vectors $(\bm{\theta},\bm{p})$ into two minimums as follows:
\begin{equation*}
\min_{(\bm{\theta},\bm{p})}
d(\bm{X}(\bm{p})\bm{\theta},\bm{y})
=
\min_{\bm{p}}
\min_{\bm{\theta}}
d(\bm{X}(\bm{p})\bm{\theta},\bm{y}).
\end{equation*}

Next, we consider the inner minimization problem
\begin{equation*}
\begin{aligned}
&&
\min_{\bm{\theta}}
&&&
d(\bm{X}(\bm{p})\bm{\theta},\bm{y}).
\end{aligned}
\end{equation*}

For any regular vector $\bm{p}$, this problem finds a best approximate solution to equation \eqref{E-Xpthetaeqy}. By applying Theorem~\ref{T-Axeqb} to this equation, we first evaluate 
\begin{equation}
\Delta(\bm{p})
=
(\bm{X}(\bm{p})(\bm{y}^{-}\bm{X}(\bm{p}))^{-})^{-}\bm{y}.
\label{E-Deltap}
\end{equation}

Furthermore, we obtain the best approximation error
\begin{equation*}
\min_{\bm{\theta}}
d(\bm{X}(\bm{p})\bm{\theta},\bm{y})
=
\sqrt{\Delta(\bm{p})},
\end{equation*}
and the solution vector
\begin{equation*}
\bm{\theta}(\bm{p})
=
\sqrt{\Delta(\bm{p})}
\left(
\bm{y}^{-}\bm{X}(\bm{p})
\right)^{-}.
\end{equation*}

To complete the solution of the initial problem at \eqref{P-mindXpthetay}, we need to find a vector of optimal exponents by solving the problem
\begin{equation}
\begin{aligned}
&&
\min_{\bm{p}}
&&&
\Delta(\bm{p}).
\end{aligned}
\label{P-minpDeltap}
\end{equation}

Suppose that a solution of the last problem is obtained to be
\begin{equation*}
\bm{p}_{\ast}
=
\arg\min_{\bm{p}}
\Delta(\bm{p})
=
\arg\min_{\bm{p}}
(\bm{X}(\bm{p})(\bm{y}^{-}\bm{X}(\bm{p}))^{-})^{-}\bm{y}.
\end{equation*}

Then, we can evaluate the optimal vector of parameters as
\begin{equation*}
\bm{\theta}_{\ast}
=
\bm{\theta}(\bm{p}_{\ast})
=
\sqrt{\Delta(\bm{p}_{\ast})}
\left(
\bm{y}^{-}\bm{X}(\bm{p}_{\ast})
\right)^{-},
\end{equation*}
which completes the construction of an approximating polynomial.

\subsection{Determination of optimal exponents}

We now investigate problem \eqref{P-minpDeltap} of minimization of the error function $\Delta(\bm{p})$ with respect to the unknown vector $\bm{p}$ of exponents. To examine the problem, we turn from the vector representation of the objective function in \eqref{P-minpDeltap} to a scalar form. We apply distributive identities for addition and minimum to rearrange the problem into a minimization problem over all partitions of the set $\{1,\ldots,M\}$ into $N$ parts each corresponding to one unknown exponent. Given a partition, the objective function to be minimized is evaluated by solving $N$ auxiliary minimization problems to find the exponents as the sum (maximum) of the minimum values in these problems. 

We start with the error function in the vector form of \eqref{E-Deltap}. To expand it into a scalar function of exponents $p_{1},\ldots,p_{N}$, we introduce the polynomial
\begin{equation*}
\varphi(p)
=
y_{1}^{-1}x_{1}^{p}\oplus\cdots\oplus y_{M}^{-1}x_{M}^{p}
=
\bigoplus_{i=1}^{M}y_{i}^{-1}x_{i}^{p}.
\end{equation*}

With this notation, we write the vector
\begin{equation*}
(\bm{y}^{-}\bm{X}(\bm{p}))^{-}
=
\left(
\begin{array}{c}
(y_{1}^{-1}x_{1}^{p_{1}}\oplus\cdots\oplus y_{M}^{-1}x_{M}^{p_{1}})^{-1}
\\
\vdots
\\
(y_{1}^{-1}x_{1}^{p_{N}}\oplus\cdots\oplus y_{M}^{-1}x_{M}^{p_{N}})^{-1}
\end{array}
\right)
=
\left(
\begin{array}{c}
(\varphi(p_{1}))^{-1}
\\
\vdots
\\
(\varphi(p_{N}))^{-1}
\end{array}
\right).
\end{equation*}

Multiplication of the vector by $\bm{X}(\bm{p})$ from the left and some algebra yield
\begin{equation*}
\bm{X}(\bm{p})(\bm{y}^{-}\bm{X}(\bm{p}))^{-}
=
\left(
\begin{array}{c}
\displaystyle{\bigoplus_{j=1}^{N}}
(x_{1}^{-p_{j}}\varphi(p_{j}))^{-1}
\\
\vdots
\\
\displaystyle{\bigoplus_{j=1}^{N}}
(x_{M}^{-p_{j}}\varphi(p_{j}))^{-1}
\end{array}
\right).
\end{equation*}

Further calculation leads to the result
\begin{equation*}
\Delta(\bm{p})
=
(\bm{X}(\bm{p})(\bm{y}^{-}\bm{X}(\bm{p}))^{-})^{-}\bm{y}
=
\bigoplus_{i=1}^{M}
\left(
\bigoplus_{j=1}^{N}
(y_{i}x_{i}^{-p_{j}}\varphi(p_{j}))^{-1}
\right)^{-1}.
\end{equation*}

After introducing the notation
\begin{equation*}
\varphi_{i}(p)
=
y_{i}x_{i}^{-p}\varphi(p),
\qquad
i=1,\ldots,M;
\end{equation*}
we arrive at the representation
\begin{equation*}
\Delta(\bm{p})
=
\bigoplus_{i=1}^{M}
\left(
\bigoplus_{j=1}^{N}
(\varphi_{i}(p_{j}))^{-1}
\right)^{-1}.
\end{equation*}

Furthermore, we observe that
\begin{equation*}
\left(
\bigoplus_{j=1}^{N}
(\varphi_{i}(p_{j}))^{-1}
\right)^{-1}
=
\min_{1\leq j\leq N}
\varphi_{i}(p_{j}),
\qquad
i=1,\ldots,M;
\end{equation*}
and thus represent the above function in the form
\begin{equation*}
\Delta(\bm{p})
=
\bigoplus_{i=1}^{M}
\min_{1\leq j\leq N}
\varphi_{i}(p_{j})
=
\bigoplus_{i=1}^{M}
\min(\varphi_{i}(p_{1}),\ldots,\varphi_{i}(p_{N})).
\end{equation*}

By applying identity \eqref{E-plusmin}, we replace this function as follows:
\begin{equation*}
\bigoplus_{i=1}^{M}
\min_{1\leq j\leq N}
\varphi_{i}(p_{j})
=
\min_{(I_{1},\ldots,I_{N})}
\bigoplus_{j=1}^{N}
\bigoplus_{i\in I_{j}}
\varphi_{i}(p_{j}),
\end{equation*}
where the minimum on the right-hand side is taken over all ordered set partitions $(I_{1},\ldots,I_{N})$ of the set of naturals $\{1,\ldots,M\}$ into $N$ parts.

We note that in the context of finding optimal exponents $p_{1},\ldots,p_{N}$ for the approximating polynomial, the order of the exponents obtained does not matter. Indeed, since the optimal value of the corresponding parameters (coefficients) $\theta_{1},\ldots,\theta_{N}$ are fixed after determination of exponents, the numbering of monomials (and hence exponents) can be changed without affecting the polynomial. As a result, we can replace the minimization over ordered set partitions $(I_{1},\ldots,I_{N})$ by minimizing over ordinary partitions $\{I_{1},\ldots,I_{N}\}$.

In accordance with the above observation, problem \eqref{P-minpDeltap} becomes
\begin{equation}
\begin{aligned}
&&
\min_{p_{1},\ldots,p_{N}}
&&&
\min_{\{I_{1},\ldots,I_{N}\}}
\bigoplus_{j=1}^{N}
\bigoplus_{i\in I_{j}}
\varphi_{i}(p_{j}),
\end{aligned}
\label{P-minp1pN}
\end{equation}
where the inner minimum is over partitions $\{I_{1},\ldots,I_{N}\}$ of the set $\{1,\ldots,M\}$.

To find the optimal vector $\bm{p}_{\ast}=(p_{1}^{\ast},\ldots,p_{N}^{\ast})^{T}$ which is the minimizer
\begin{equation*}
(p_{1}^{\ast},\ldots,p_{N}^{\ast})
=
\arg\min_{p_{1},\ldots,p_{N}}
\min_{\{I_{1},\ldots,I_{N}\}}
\bigoplus_{j=1}^{N}
\bigoplus_{i\in I_{j}}
\varphi_{i}(p_{j}),
\end{equation*}
we change the order of minimization, and arrive at the equivalent problem 
\begin{equation*}
\begin{aligned}
&&
\min_{\{I_{1},\ldots,I_{N}\}}
&&&
\min_{p_{1},\ldots,p_{N}}
\bigoplus_{j=1}^{N}
\bigoplus_{i\in I_{j}}
\varphi_{i}(p_{j}).
\end{aligned}
\end{equation*}

Finally, application of identity \eqref{E-minplus} yields the problem
\begin{equation}
\begin{aligned}
&&
\min_{\{I_{1},\ldots,I_{N}\}}
&&&
\bigoplus_{j=1}^{N}
\min_{p_{j}}
\bigoplus_{i\in I_{j}}
\varphi_{i}(p_{j}),
\end{aligned}
\label{P-minI1INminpjvarphiipj}
\end{equation}
where the function under summation is given by
\begin{equation*}
\varphi_{i}(p_{j})
=
y_{i}x_{i}^{-p_{j}}
(y_{1}^{-1}x_{1}^{p_{j}}\oplus\cdots\oplus y_{M}^{-1}x_{M}^{p_{j}}).
\end{equation*}

\subsection{Representation of solution}

Consider the polynomial approximation problem under study, which is now represented as \eqref{P-minI1INminpjvarphiipj}. Denote the minimum of the objective function in the problem by $\Delta_{\ast}$ and the partition that attains this minimum by $\{I_{1}^{\ast},\ldots,I_{N}^{\ast}\}$.

Then, the exponents of the best approximating polynomial are defined as
\begin{equation*}
p_{1}^{\ast}
=
\arg\min_{p}
\bigoplus_{i\in I_{1}^{\ast}}
\varphi_{i}(p),
\quad
\ldots,
\quad
p_{N}^{\ast}
=
\arg\min_{p}
\bigoplus_{i\in I_{N}^{\ast}}
\varphi_{i}(p),
\end{equation*}
whereas the corresponding parameters of the polynomial are
\begin{equation*}
\theta_{1}^{\ast}
=
\sqrt{\Delta_{\ast}}(\varphi(p_{1}^{\ast}))^{-1},
\quad
\ldots,
\quad
\theta_{N}^{\ast}
=
\sqrt{\Delta_{\ast}}(\varphi(p_{N}^{\ast}))^{-1}.
\end{equation*}

We observe that the development of a workable solution to problem \eqref{P-minI1INminpjvarphiipj} requires to take into account several issues. Specifically, an efficient method of solving the inner minimization problems at \eqref{P-minI1INminpjvarphiipj}, which is the key component of the solution, needs to be implemented. Another issue is combinatorial number of partitions to examine in problem \eqref{P-minI1INminpjvarphiipj}, which may make the problem hard to solve for large $M$ and $N$. To address these issues, we further transform the problem and propose an economical procedure to examine partitions in the next section.

\section{Approximation by max-plus and max-algebra polynomials}
\label{S-AMPMAP}

Since the solution of the approximation problem at \eqref{P-minI1INminpjvarphiipj} in the setting of an arbitrary idempotent semifield may be difficult, we concentrate in this section on the solution in terms of max-plus algebra and max-algebra. We examine the inner minimization problem in \eqref{P-minI1INminpjvarphiipj} to develop a reasonable procedure of finding this minimum.

For a fixed partition $\{I_{1},\ldots,I_{N}\}$ of the set $\{1,\ldots,M\}$, consider a problem that is given by a nonempty part $I$ of the partition and written as 
\begin{equation}
\begin{aligned}
&&
\min_{p}
&&&
\bigoplus_{i\in I}
\varphi_{i}(p).
\end{aligned}
\label{P-minpvarphiip}
\end{equation}

Below we show how this problem can be rearranged and solved in the framework of max-plus algebra, and then extend the solution to max-algebra.

\subsection{Approximation by max-plus polynomials}

We start with the development of solution of the approximation problem at \eqref{P-minI1INminpjvarphiipj} in the context of the semifield $\mathbb{R}_{\max,+}$ (max-plus algebra), where addition $\oplus$ is defined as $\max$ and multiplication $\otimes$ as $+$.

Let us verify that in this case the objective function in problem \eqref{P-minpvarphiip} can be transformed into a polynomial in which $p$ is its indeterminate. Indeed, in terms of max-plus algebra, we can write the equality $x^{y}=y^{x}$ for all $x,y\in\mathbb{R}$ since both sides of this equality correspond to the same product $x\times y$ in the ordinary notation. Therefore, we have two equivalent representations for the functions
\begin{gather*}
\varphi(p)
=
y_{1}^{-1}x_{1}^{p}\oplus\cdots\oplus y_{M}^{-1}x_{M}^{p}
=
y_{1}^{-1}p^{x_{1}}\oplus\cdots\oplus y_{M}^{-1}p^{x_{M}},
\\
\varphi_{i}(p)
=
y_{i}x_{i}^{-p}
(y_{1}^{-1}x_{1}^{p}\oplus\cdots\oplus y_{M}^{-1}x_{M}^{p})
=
y_{i}p^{-x_{i}}
(y_{1}^{-1}p^{x_{1}}\oplus\cdots\oplus y_{M}^{-1}p^{x_{M}});
\end{gather*}
where $p$ first plays the role of an exponent, and then of an indeterminate.

As a result, we arrive at polynomials in $p$, which take the form
\begin{gather*}
\varphi(p)
=
\bigoplus_{i=1}^{M}y_{i}^{-1}p^{x_{i}},
\\
\varphi_{i}(p)
=
y_{1}^{-1}y_{i}p^{x_{1}-x_{i}}\oplus\cdots\oplus y_{M}^{-1}y_{i}p^{x_{M-x_{i}}}
=
\bigoplus_{j=1}^{M}
y_{j}^{-1}y_{i}
p^{x_{j}-x_{i}}.
\end{gather*}

Consider the objective function in problem \eqref{P-minpvarphiip}, given by
\begin{equation*}
\bigoplus_{i\in I}
\varphi_{i}(p)
=
\bigoplus_{i\in I}
\bigoplus_{j=1}^{M}
y_{j}^{-1}y_{i}
p^{x_{j}-x_{i}},
\end{equation*}
which is a sum of polynomials in $p$ and hence a polynomial itself.

We rearrange this polynomial to reorder the monomials in increasing exponents and merge the monomials with equal exponents. We renumber the monomials and denote by $u_{j}$ the coefficient and by $v_{j}$ the exponent of monomial $j$ for each $j=1,\ldots,L$, where $L\leq M|I|$ is the number of monomials in the obtained polynomial. The objective function becomes
\begin{equation*}
\bigoplus_{i\in I}
\varphi_{i}(p)
=
\bigoplus_{j=1}^{L}
u_{j}p^{v_{j}}.
\end{equation*}

Finally, problem \eqref{P-minpvarphiip} takes the form of a polynomial optimization problem, which is given by
\begin{equation*}
\begin{aligned}
&&
\min_{p}
&&&
\bigoplus_{j=1}^{L}
u_{j}p^{v_{j}}.
\end{aligned}
\end{equation*}

As one can see, the minimum value of the objective function and the solution of the problem can be found by application of Lemma~\ref{L-minxajxpj}.

We turn back to problem \eqref{P-minI1INminpjvarphiipj} and note that this problem involves the solution of a set of polynomial minimization problems in one variable. At the same time, direct solution of all these problems leads to a combinatorial growth of the polynomials to handle. Indeed, the number of the minimization problems to solve is determined by the number of partitions $\{I_{1},\ldots,I_{N}\}$, which is given by a Stirling number, whereas each partition leads to the solution of one to $N$ problems. To overcome this difficulty, we apply a solution scheme that follows the key ideas of agglomerative (hierarchical) clustering \cite{Ward1963Hierarchical,Anderberg1973Cluster}.

\subsection{Derivation of optimal exponents}

We propose a procedure that starts with a partition $\{I_{1},\ldots,I_{M}\}$ of $M\geq N$ one-element subsets $I_{1}=\{1\},\ldots,I_{M}=\{M\}$ and then gradually reduces the partition by merging the subsets in pairs until the number of subsets becomes $N$. In the context of the optimization problem under consideration, each subset in a partition determines a polynomial in one variable to be minimized. Specifically, the initial partition $\{I_{1},\ldots,I_{M}\}$ corresponds to the solution of $M$ polynomial optimization problems with the objective functions $\varphi_{j}(p_{j})$ in variables $p_{j}$ for all $j=1,\ldots,M$. The minimum value and solution set of the optimization problem incident to a subset become the characteristics of the subset that are utilized by the procedure to make further steps of solution. The entire partition is then characterized by the maximum of the minimums for the partition subsets, whereas the solutions for the subsets determine related components of the solution vector.

If two partition subsets are merged into a new subset, the corresponding polynomials are combined (added) to form a single polynomial (which is further simplified by merging monomials with a common exponent if any occur). Since the polynomials are monotone functions, the new polynomial may have a minimum that is greater or equal to the maximum of minimums of the initial polynomials, as well as may have a new solution set.

At each step of the procedure, those two subsets in the current partition are selected for merging whose merged polynomial has the least minimum value over all pairs of the subsets. As a result, the objective function in problem \eqref{P-minI1INminpjvarphiipj}, which is evaluated as the maximum of the minimums of involved polynomials, may increase its value or does not change. This offers a reasonable way to solve the problem by reducing the number of partitions according to the minimum increase of the objective function per merge. Upon completion of the procedure, there are $N$ polynomials whose minimum values are maximized to evaluate the approximation error. The set of solutions of minimization problems for all polynomials is taken as a set of optimal values of the exponents $p_{1},\ldots,p_{N}$.   

In terms of agglomerative clustering \cite{Ward1963Hierarchical,Anderberg1973Cluster}, the initial objects to be grouped into $N$ clusters are the polynomials $\varphi_{j}$ for all $j=1,\ldots,M$. The merge of two clusters corresponds to the formation of a new polynomial that is given by the sum of the polynomials associated with the merged clusters. The similarity measure (distance function) between two clusters is defined as the minimum value of a polynomial that results from merging these clusters. At each step, the procedure combines two clusters that are most similar into one bigger cluster, which results in a new polynomial. The procedure ends when the number of clusters (polynomials) decreases to $N$.

\subsection{Approximation procedure}

The approximation procedure can be described in the form of Algorithm~\ref{A-AMPP} based on the partition of the set $\{1,\ldots,M\}$ by using agglomerative clustering. The main component of the algorithm is the solution of polynomial minimization problems at Step~\ref{Loop} by using the result of Lemma~\ref{L-minxajxpj}.
\begin{algorithm}
\caption{Approximation by max-plus polynomial}
\label{A-AMPP}
\bigskip
\begin{description}
\item[Step 0:]
Input $M$ samples $(x_{i},y_{i})$ for $i=1,\ldots,M$ and fix a positive integer $N$.
Define the max-plus algebra polynomials
\begin{equation*}
\varphi_{i}(p)
=
\bigoplus_{j=1}^{M}
y_{j}^{-1}y_{i}
p^{x_{j}-x_{i}},
\qquad
i=1,\ldots,M.
\end{equation*}
Construct the subsets and form the partition
\begin{equation*}
I_{1}^{(0)}
=
\{1\},
\ldots,
I_{M}^{(0)}
=
\{M\};
\qquad
\mathcal{P}^{(0)}
=
\{I_{1}^{(0)},\ldots,I_{M}^{(0)}\};
\end{equation*}
and put $k=1$.
\item[Step 1:]\label{Loop}
Find a pair of subsets to satisfy the condition
\begin{equation*}
(U^{(k)},V^{(k)})
=
\arg\min_{\substack{U,V\in\mathcal{P}^{(k-1)}\\ U\ne V}}
\min_{p}
\left(
\bigoplus_{i\in U}
\varphi_{i}(p)
\oplus
\bigoplus_{i\in V}
\varphi_{i}(p)
\right).
\end{equation*}
\item[Step 2:]
Merge the subsets $U^{(k)}$ and $V^{(k)}$ to form the partition 
\begin{equation*}
\mathcal{P}^{(k)}
=
(\mathcal{P}^{(k-1)}\setminus\{U^{(k)},V^{(k)}\})\cup\{U^{(k)}\cup V^{(k)}\}.
\end{equation*}
\item[Step 3:]
If $M-k=N$, then define an optimal partition
\begin{equation*}
\mathcal{P}_{\ast}
=
\{I_{1}^{\ast},\ldots,I_{N}^{\ast}\}
=
\mathcal{P}^{(M-k)}
=
\{I_{1}^{(M-k)},\ldots,I_{N}^{(M-k)}\};
\end{equation*}
otherwise put $k=k+1$ and go to Step~\ref{Loop}.
\item[Step 4:]
Derive the minimums and exponents corresponding to $\mathcal{P}_{\ast}$ to be
\begin{gather*}
\delta_{1}^{\ast}
=
\min_{p}
\bigoplus_{i\in I_{1}^{\ast}}
\varphi_{i}(p),
\quad
\ldots,
\quad
\delta_{N}^{\ast}
=
\min_{p}
\bigoplus_{i\in I_{N}^{\ast}}
\varphi_{i}(p);
\\
p_{1}^{\ast}
=
\arg\min_{p}
\bigoplus_{i\in I_{1}^{\ast}}
\varphi_{i}(p),
\quad
\ldots,
\quad
p_{N}^{\ast}
=
\arg\min_{p}
\bigoplus_{i\in I_{N}^{\ast}}
\varphi_{i}(p).
\end{gather*}
Evaluate the squared approximation error
\begin{equation*}
\Delta_{\ast}
=
\bigoplus_{j=1}^{N}
\delta_{j}^{\ast}
=
\max_{1\leq j\leq N}
\delta_{j}^{\ast},
\end{equation*}
and calculate the parameters 
\begin{equation*}
\theta_{1}^{\ast}
=
\sqrt{\Delta_{\ast}}(\varphi(p_{1}^{\ast}))^{-1},
\quad
\ldots,
\quad
\theta_{N}^{\ast}
=
\sqrt{\Delta_{\ast}}(\varphi(p_{N}^{\ast}))^{-1}.
\end{equation*}
\end{description}
\end{algorithm}

To estimate the computational complexity of Algorithm~\ref{A-AMPP}, we derive an upper bound of the time required for calculation. We note that the most computationally intensive part of the procedure is the solution of optimization problems at Step~\ref{Loop}.

At the first iteration with $k=1$, the algorithm solves $M(M-1)/2$ initial minimization problems for all pairwise sums of polynomials. Each sum to minimize is a polynomial in one variable with no more than $2M$ monomials, which can be minimized by applying Lemma~\ref{L-minxajxpj} in no more than $O(M^{2})$ operations. Therefore, the computational complexity of the first iteration can be estimated as $O(M^{4})$.

As a result of this iteration, two polynomials are selected for which minimizing the sum of them yields the minimum value over all pairs of polynomials. The selected polynomials are replaced by their sum as a new polynomial in the set, which decrements the overall number of polynomials by one. Furthermore, at iteration $k=2$, the algorithm needs to solve $M-2$ minimization problems for the sums of this new polynomial with other $M-2$ polynomials that remain from the previous iteration together with the results of the minimization of their pairwise sums.

Consider iteration $k>1$, and note that the overall number of iterations of the algorithm is not greater than $M-N$. The number of polynomial minimization problems to solve at iteration $k$ is $M-k-1<M$. The number of monomials in each polynomial to minimize does not exceed $(M-N)M$, and thus minimization of the polynomial involves no more than $O((M-N)^{2}M^{2})$ operations. By combining these estimates, we conclude that the time required to perform the iterations $k=2,\ldots,N-M$ is no more than $O((M-N)^{3}M^{3})$. As a result, the computational complexity of the entire algorithms can be estimated as $O(\max(M^{4},(M-N)^{3}M^{3}))$.

Finally, we note that the agglomerative clustering scheme is known to yield greedy algorithms that do not guarantee a globally optimal solution. However, as numerical examples in the next section show, the application of this clustering technique in the approximation problems under consideration is quite capable of providing results with a reasonable accuracy.

\subsection{Approximation by max-algebra polynomials}

Suppose that problem \eqref{P-minI1INminpjvarphiipj} is solved in the framework of the tropical semifield $\mathbb{R}_{\max}$ (max-algebra), where the tropical addition $\oplus$ is defined as $\max$, and multiplication $\otimes$ as $\times$. We note that, in contrast to max-plus algebra, the tropical powers $x^{y}$ and $y^{x}$ now have different values and cannot interchange, which does not allow the same transformation of exponents as before.

Due to the monotonicity of logarithm, the solution of \eqref{P-minI1INminpjvarphiipj} is equivalent to the minimization over all partitions $\{I_{1},\ldots,I_{N}\}$ of the function
\begin{equation*}
\log
\bigoplus_{j=1}^{N}
\min_{p_{j}}
\bigoplus_{i\in I_{j}}
\varphi_{i}(p_{j})
=
\bigoplus_{j=1}^{N}
\min_{p_{j}}
\bigoplus_{i\in I_{j}}
\log\varphi_{i}(p_{j}),
\end{equation*}
where the logarithm to the natural base is assumed.

Since the function $\log\varphi_{i}(p)$ is a max-plus polynomial, we can reduce the problem of approximation by max-polynomials to approximation by max-plus polynomials described in the previous section. The solution includes additional steps for pre- and post-processing of data. The pre-processing consists in the replacement of the initial samples $x_{1},\ldots,x_{M}$ and $y_{1},\ldots,y_{M}$ by $\log x_{1},\ldots,\log x_{M}$ and $\log y_{1},\ldots,\log y_{M}$. After the completion of the algorithm, we replace the obtained approximation error $\Delta_{\ast}$ by $\exp(\Delta_{\ast})$ and the parameters $\theta_{1}^{\ast},\ldots,\theta_{N}^{\ast}$ by $\exp(\theta_{1}^{\ast}),\ldots,\exp(\theta_{N}^{\ast})$.

We observe that the pre-processing step involves $2M$ additional calculations of logarithmic function, whereas post-processing does $N+1$ calculations of exponential function. Since the number of operations required for the pre- and post-processing steps has a linear growth with respect to $M$ and $N$, these steps do not sufficiently increase the time complexity of Algorithm~\ref{A-AMPP}. As a result, we conclude that both solutions of max-plus and max-algebra polynomial approximation problems have the same rate of complexity.

\section{Illustrative numerical examples}
\label{S-INE}

In this section, we illustrate the obtained results with numerical examples of approximation problems using max-plus and max-algebra polynomials. As a test problem, we consider the approximation of the function defined in terms of the standard algebra as
\begin{equation*}
f(x)
=
(x-3/4)^{2}-3(x-1)^{1/2}+2,
\qquad
x\in[1,3].
\end{equation*}

Given $M=21$ samples $x_{i}=1+(i-1)/10$ and $y_{i}=f(x_{i})$ of the input and output of the above function for $i=1,\ldots,M$, the problem is to approximate this function by means of polynomials defined in the max-plus and max-algebra settings. Below we demonstrate results for polynomials with different number of monomials.

We start with the solutions based on the implementation of Puiseux polynomials in the max-plus algebra setting. As a result of approximation by polynomials with $N=5$ monomials, Algorithm~\ref{A-AMPP} yields the squared approximation error $\Delta_{\ast}=0.1054$. The vectors of exponents $\bm{p}_{\ast}=(p_{1}^{\ast},\ldots,p_{N}^{\ast})^{T}$ and coefficients $\bm{\theta}_{\ast}=(\theta_{1}^{\ast},\ldots,\theta_{N}^{\ast})^{T}$ obtained are as follows:
\begin{gather*}
\bm{p}_{\ast}
=
\left(
\begin{array}{ccccc}
-8.8868 & -1.3989 & 0.2591 & 1.3982 & 2.7450
\end{array}
\right)^{T},
\\
\bm{\theta}_{\ast}
=
\left(
\begin{array}{ccccc}
11.0020 & 2.5306 & 0.0047 & -2.2076 & -5.4679
\end{array}
\right)^{T}.
\end{gather*}

The approximating function is written using standard arithmetic operations as
\begin{multline*}
P_{\ast}(x)
=
\max(
-8.8868x+11.0020,
-1.3989x+2.5306,
\\
0.2591x+0.0047,
1.3982x-2.2076,
2.7450x-5.4679).
\end{multline*}

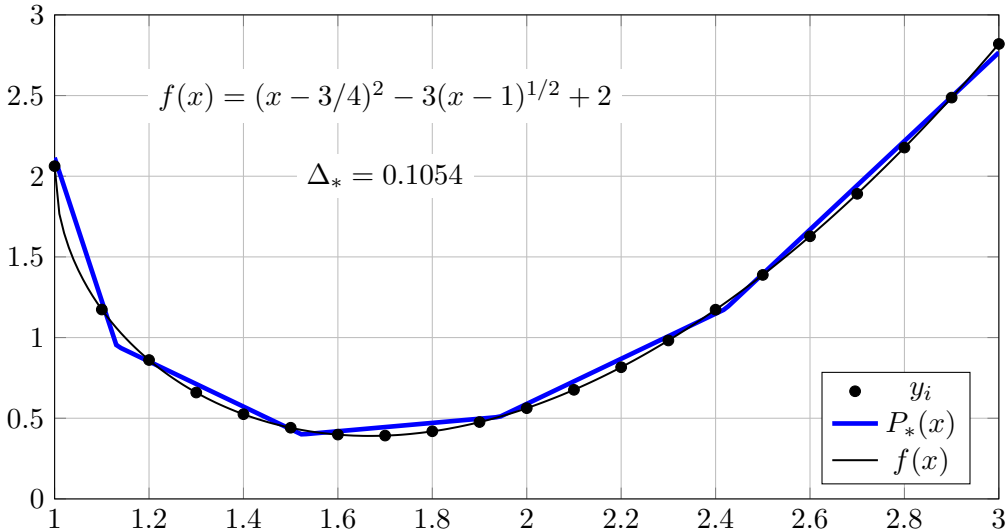
\begin{figure}[ht]
\begin{tikzpicture}

\begin{axis}[legend pos=south east,
width=14cm,
height=8cm,
grid=major,
ymin=0,
ymax=3,
xmin=1,
xmax=3
]

\addplot[
only marks,
]
coordinates {
(1.0000,2.0625)
(1.1000,1.1738)
(1.2000,0.8609)
(1.3000,0.6593)
(1.4000,0.5251)
(1.5000,0.4412)
(1.6000,0.3987)
(1.7000,0.3925)
(1.8000,0.4192)
(1.9000,0.4765)
(2.0000,0.5625)
(2.1000,0.6761)
(2.2000,0.8162)
(2.3000,0.9820)
(2.4000,1.1729)
(2.5000,1.3883)
(2.6000,1.6278)
(2.7000,1.8910)
(2.8000,2.1776)
(2.9000,2.4873)
(3.0000,2.8199)
};

\addlegendentry{$y_{i}$}

\addplot[
blue,
samples=200,
line width=1.75pt,
domain=1.0:3.0,
y domain=0:10
]{max(
11.0020+x*(-8.8868),
 2.5306+x*(-1.3989),
 0.0047+x*0.2591,   
-2.2076+x*1.3982,   
-5.4679+x*2.7450   
)};

\addlegendentry{$P_{\ast}(x)$}

\addplot[
samples=200,
black,
line width=0.75pt,
domain=1.0:3.0,
y domain=0:10
]{(x-3/4)^2-3*(x-1)^(1/2)+2};

\addlegendentry{$f(x)$}

\node[style={fill=white}] at (axis cs: 1.7,2.5) {$f(x)=(x-3/4)^{2}-3(x-1)^{1/2}+2$};
\node[style={fill=white}] at (axis cs: 1.7,2.0) {$\Delta_{\ast}=0.1054$};

\end{axis}

\end{tikzpicture}
\caption{Approximation by max-plus polynomial $P_{\ast}(x)$ with $N=5$ terms.}
\label{F-AFMPPN5}
\end{figure}

A graphical illustration of the solution is given in Figure~\ref{F-AFMPPN5}.

Furthermore, we consider the result for polynomials with $N=7$ monomials. In this case, the squared error is $\Delta_{\ast}=0.0321$, whereas the vectors of exponents and coefficients are
\begin{gather*}
\bm{p}_{\ast}
=
\left(
\begin{array}{ccccccc}
-8.8868 & -2.0153 & -0.8395 & 0.2591 & 1.3982 & 2.3938 & 3.2114
\end{array}
\right)^{T},
\\
\bm{\theta}_{\ast}
=
\left(
\begin{array}{ccccccc}
10.9654 & 3.2952 & 1.7165 & -0.0320 & -2.2442 & -4.5801 & -6.8097
\end{array}
\right)^{T}.
\end{gather*}

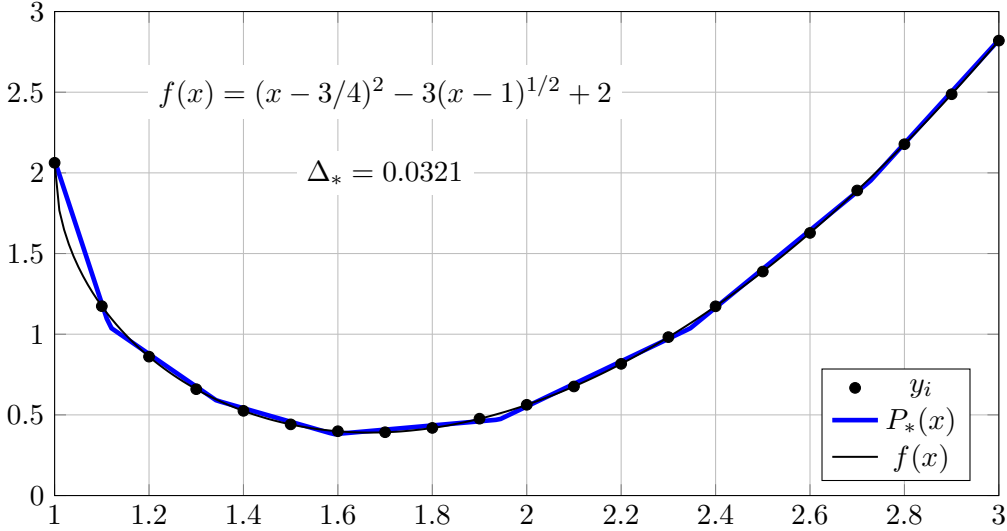
\begin{figure}[ht]
\begin{tikzpicture}

\begin{axis}[legend pos=south east,
width=14cm,
height=8cm,
grid=major,
ymin=0,
ymax=3,
xmin=1,
xmax=3
]

\addplot[
only marks,
]
coordinates {
(1.0000,2.0625)
(1.1000,1.1738)
(1.2000,0.8609)
(1.3000,0.6593)
(1.4000,0.5251)
(1.5000,0.4412)
(1.6000,0.3987)
(1.7000,0.3925)
(1.8000,0.4192)
(1.9000,0.4765)
(2.0000,0.5625)
(2.1000,0.6761)
(2.2000,0.8162)
(2.3000,0.9820)
(2.4000,1.1729)
(2.5000,1.3883)
(2.6000,1.6278)
(2.7000,1.8910)
(2.8000,2.1776)
(2.9000,2.4873)
(3.0000,2.8199)
};

\addlegendentry{$y_{i}$}

\addplot[
blue,
samples=200,
line width=1.75pt,
domain=1.0:3.0,
y domain=0:10
]{max(
10.9654+x*(-8.8868),
 3.2952+x*(-2.0153),
 1.7165+x*(-0.8395),
-0.0320+x*0.2591,  
-2.2442+x*1.3982,  
-4.5801+x*2.3938,  
-6.8097+x*3.2114   
)};

\addlegendentry{$P_{\ast}(x)$}

\addplot[
samples=200,
black,
line width=0.75pt,
domain=1.0:3.0,
y domain=0:10
]{(x-3/4)^2-3*(x-1)^(1/2)+2};

\addlegendentry{$f(x)$}

\node[style={fill=white}] at (axis cs: 1.7,2.5) {$f(x)=(x-3/4)^{2}-3(x-1)^{1/2}+2$};
\node[style={fill=white}] at (axis cs: 1.7,2.0) {$\Delta_{\ast}=0.0321$};

\end{axis}

\end{tikzpicture}
\caption{Approximation by max-plus polynomial $P_{\ast}(x)$ with $N=7$ terms.}
\label{F-AFMPPN7}
\end{figure}

The results obtained are shown Figure~\ref{F-AFMPPN7}, where the best approximate function can be represented in terms of standard operations as
\begin{multline*}
P_{\ast}(x)
=
\max(
-8.8868x+10.9654,
-2.0153x+3.2952,
-0.8395x+1.7165,
\\
0.2591x-0.0320,
1.3982x-2.2442,
2.3938x-4.5801,
3.2114x-6.8097).
\end{multline*}

\begin{figure}[ht]
\begin{tikzpicture}

\begin{axis}[legend pos=south east,
width=14cm,
height=8cm,
grid=major,
ymin=0,
ymax=3,
xmin=1,
xmax=3
]

\addplot[
only marks,
]
coordinates {
(1.0000,2.0625)
(1.1000,1.1738)
(1.2000,0.8609)
(1.3000,0.6593)
(1.4000,0.5251)
(1.5000,0.4412)
(1.6000,0.3987)
(1.7000,0.3925)
(1.8000,0.4192)
(1.9000,0.4765)
(2.0000,0.5625)
(2.1000,0.6761)
(2.2000,0.8162)
(2.3000,0.9820)
(2.4000,1.1729)
(2.5000,1.3883)
(2.6000,1.6278)
(2.7000,1.8910)
(2.8000,2.1776)
(2.9000,2.4873)
(3.0000,2.8199)
};

\addlegendentry{$y_{i}$}

\addplot[
blue,
samples=200,
line width=1.75pt,
domain=1.0:3.0,
y domain=0:10
]{max(
10.9493+x*(-8.8868),
 3.2792+x*(-2.0153),
 1.7005+x*(-0.8395),
 0.4978+x*(-0.0619),
-0.6110+x*0.5723,   
-1.7090+x*1.1357,   
-2.8316+x*1.6581,   
-3.9971+x*2.1541,   
-5.2157+x*2.6321,   
-6.4942+x*3.0971,   
-7.1574+x*3.3257   
)};

\addlegendentry{$P_{\ast}(x)$}

\addplot[
samples=200,
black,
line width=0.75pt,
domain=1.0:3.0,
y domain=0:10
]{(x-3/4)^2-3*(x-1)^(1/2)+2};

\addlegendentry{$f(x)$}

\node[style={fill=white}] at (axis cs: 1.7,2.5) {$f(x)=(x-3/4)^{2}-3(x-1)^{1/2}+2$};
\node[style={fill=white}] at (axis cs: 1.7,2.0) {$\Delta_{\ast}=0$};

\end{axis}

\end{tikzpicture}
\caption{Approximation by max-plus polynomial $P_{\ast}(x)$ with $N=11$ terms.}
\label{F-AFMPPN11}
\end{figure}
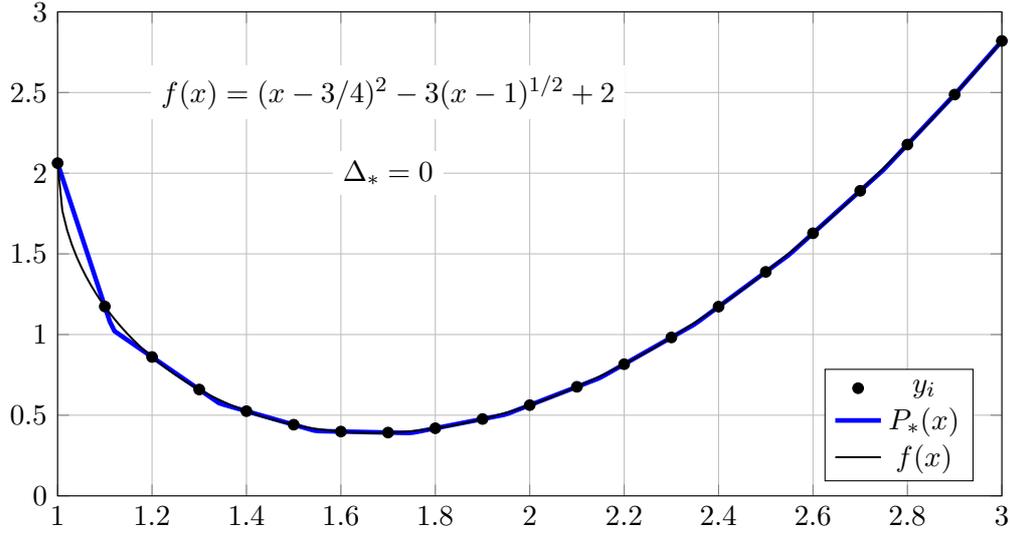

Finally, we present an approximate polynomial that consists of $N=11$ monomials (see Figure~\ref{F-AFMPPN11}). Application of Algorithm~\ref{A-AMPP} leads to the squared error $\Delta_{\ast}=0$ and the approximating function represented in the conventional form as
\begin{multline*}
P_{\ast}(x)
=
\max(
-8.8868x+10.9493,
-2.0153x+3.2792,
-0.8395x+1.7005,
\\
-0.0619x+0.4978,
0.5723x-0.6110, 
1.1357x-1.7090, 
1.6581x-2.8316, 
\\
2.1541x-3.9971, 
2.6321x-5.2157, 
3.0971x-6.4942, 
3.3257x-7.1574).
\end{multline*}

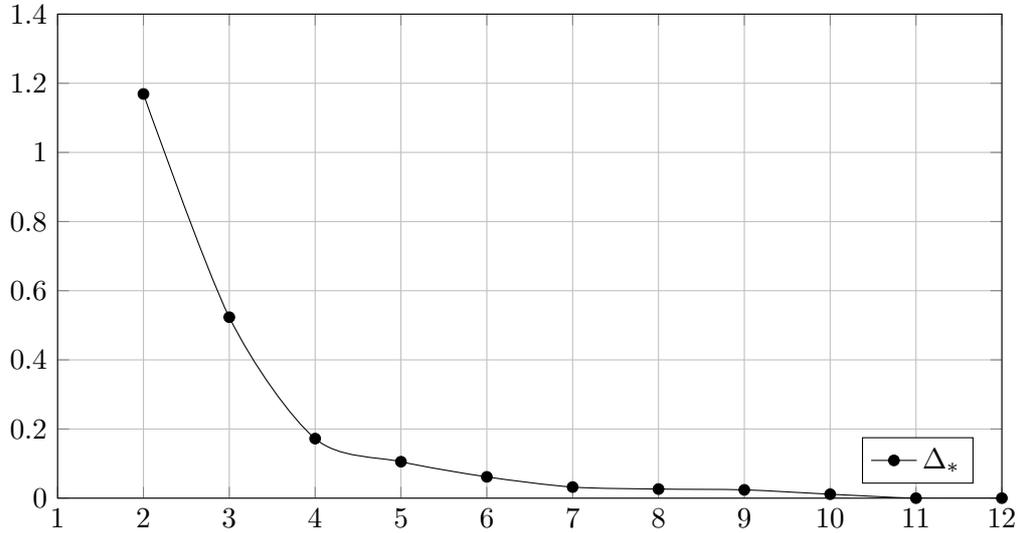
\begin{figure}[ht]
\begin{tikzpicture}

\begin{axis}[legend pos=south east,
width=14cm,
height=8cm,
grid=major,
ymin=0,
ymax=1.4,
xmin=1,
xmax=12
]

\addplot[
mark=*,
smooth
]
coordinates {
(2,1.1690)
(3,0.5234)
(4,0.1722)
(5,0.1054)
(6,0.0616)
(7,0.0321)
(8,0.0263)
(9,0.0240)
(10,0.0114)
(11,0)
(12,0)
};

\addlegendentry{$\Delta_{\ast}$}


\end{axis}

\end{tikzpicture}
\caption{Approximation error with respect to the number $N$ of monomials.}
\label{F-AEWRNM}
\end{figure}

As one can expect, the accuracy of approximation improves when the number of monomials in the approximating polynomials increases. Figure~\ref{F-AEWRNM} shows the dynamics of the squared approximation error with respect to the number $N$ of monomials in the approximating polynomial. 

We now consider the implementation of Algorithm~\ref{A-AMPP} to the approximation by Puiseux polynomials in terms of max-algebra. The approximation by polynomials with $N=5$ results in the squared approximation error $\Delta_{\ast}=1.0330$ and the vectors
\begin{gather*}
\bm{p}_{\ast}
=
\left(
\begin{array}{ccccc}
-5.9139 & -2.9957 & -0.2581 & 2.3669 & 4.0040
\end{array}
\right)^{T},
\\
\bm{\theta}_{\ast}
=
\left(
\begin{array}{ccccc}
2.0963 & 1.4625 & 0.4575 & 0.1060 & 0.0352
\end{array}
\right)^{T}.
\end{gather*}

The approximating function in the standard notation takes the form
\begin{multline*}
P_{\ast}(x)
=
\max(
2.0963x^{-5.9139},
1.4625x^{-2.9957},
\\
0.4575x^{-0.2581},
0.1060x^{2.3669},  
0.0352x^{4.0040}).
\end{multline*}

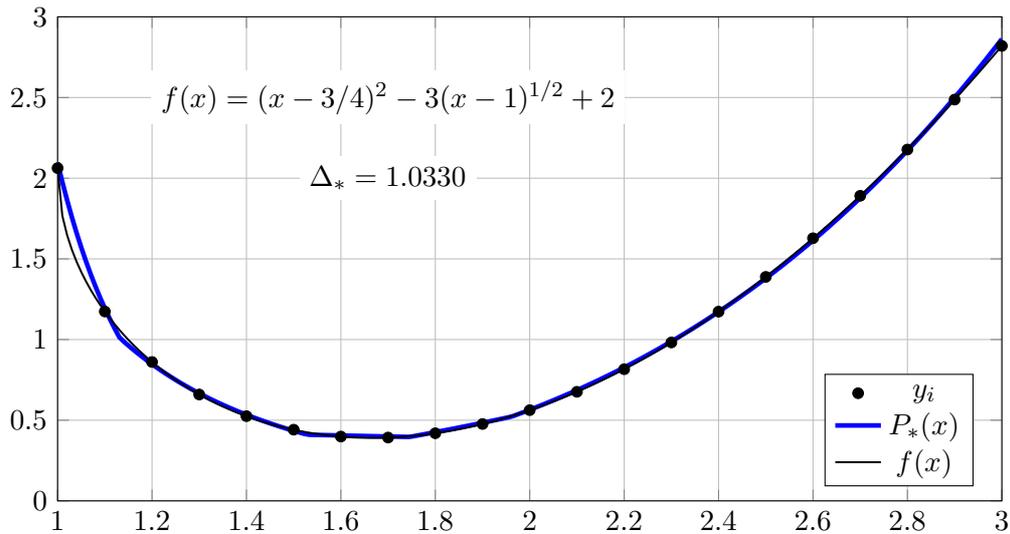
\begin{figure}[ht]
\begin{tikzpicture}

\begin{axis}[legend pos=south east,
width=14cm,
height=8cm,
grid=major,
ymin=0,
ymax=3,
xmin=1,
xmax=3
]

\addplot[
only marks,
]
coordinates {
(1.0000,2.0625)
(1.1000,1.1738)
(1.2000,0.8609)
(1.3000,0.6593)
(1.4000,0.5251)
(1.5000,0.4412)
(1.6000,0.3987)
(1.7000,0.3925)
(1.8000,0.4192)
(1.9000,0.4765)
(2.0000,0.5625)
(2.1000,0.6761)
(2.2000,0.8162)
(2.3000,0.9820)
(2.4000,1.1729)
(2.5000,1.3883)
(2.6000,1.6278)
(2.7000,1.8910)
(2.8000,2.1776)
(2.9000,2.4873)
(3.0000,2.8199)
};

\addlegendentry{$y_{i}$}

\addplot[
blue,
samples=200,
line width=1.75pt,
domain=1.0:3.0,
y domain=0:10
]{max(
2.0963*x^(-5.9139),
1.4625*x^(-2.9957),
0.4575*x^(-0.2581),
0.1060*x^2.3669,  
0.0352*x^4.0040   
)};

\addlegendentry{$P_{\ast}(x)$}

\addplot[
samples=200,
black,
line width=0.75pt,
domain=1.0:3.0,
y domain=0:10
]{(x-3/4)^2-3*(x-1)^(1/2)+2};

\addlegendentry{$f(x)$}

\node[style={fill=white}] at (axis cs: 1.7,2.5) {$f(x)=(x-3/4)^{2}-3(x-1)^{1/2}+2$};
\node[style={fill=white}] at (axis cs: 1.7,2.0) {$\Delta_{\ast}=1.0330$};

\end{axis}

\end{tikzpicture}
\caption{Approximation by max-polynomial $P_{\ast}(x)$ with $N=5$ terms.}
\label{F-AFMPN5}
\end{figure}

A graphical illustration of the solution is given in Figure~\ref{F-AFMPN5}.

\begin{figure}[ht]
\begin{tikzpicture}

\begin{axis}[legend pos=south east,
width=14cm,
height=8cm,
grid=major,
ymin=0,
ymax=3,
xmin=1,
xmax=3
]

\addplot[
only marks,
]
coordinates {
(1.0000,2.0625)
(1.1000,1.1738)
(1.2000,0.8609)
(1.3000,0.6593)
(1.4000,0.5251)
(1.5000,0.4412)
(1.6000,0.3987)
(1.7000,0.3925)
(1.8000,0.4192)
(1.9000,0.4765)
(2.0000,0.5625)
(2.1000,0.6761)
(2.2000,0.8162)
(2.3000,0.9820)
(2.4000,1.1729)
(2.5000,1.3883)
(2.6000,1.6278)
(2.7000,1.8910)
(2.8000,2.1776)
(2.9000,2.4873)
(3.0000,2.8199)
};

\addlegendentry{$y_{i}$}

\addplot[
blue,
samples=200,
line width=1.75pt,
domain=1.0:3.0,
y domain=0:10
]{max(
2.0869*x^(-5.9139),
1.5991*x^(-3.3320),
1.2426*x^(-2.5249),
0.4555*x^(-0.2581),
0.1055*x^2.3669,   
0.0351*x^4.0040,   
0.0351*x^4.0040   
)};

\addlegendentry{$P_{\ast}(x)$}

\addplot[
samples=200,
black,
line width=0.75pt,
domain=1.0:3.0,
y domain=0:10
]{(x-3/4)^2-3*(x-1)^(1/2)+2};

\addlegendentry{$f(x)$}

\node[style={fill=white}] at (axis cs: 1.7,2.5) {$f(x)=(x-3/4)^{2}-3(x-1)^{1/2}+2$};
\node[style={fill=white}] at (axis cs: 1.7,2.0) {$\Delta_{\ast}=1.0238$};

\end{axis}

\end{tikzpicture}
\caption{Approximation by max-polynomial $P_{\ast}(x)$ with $N=7$ terms.}
\label{F-AFMPN7}
\end{figure}
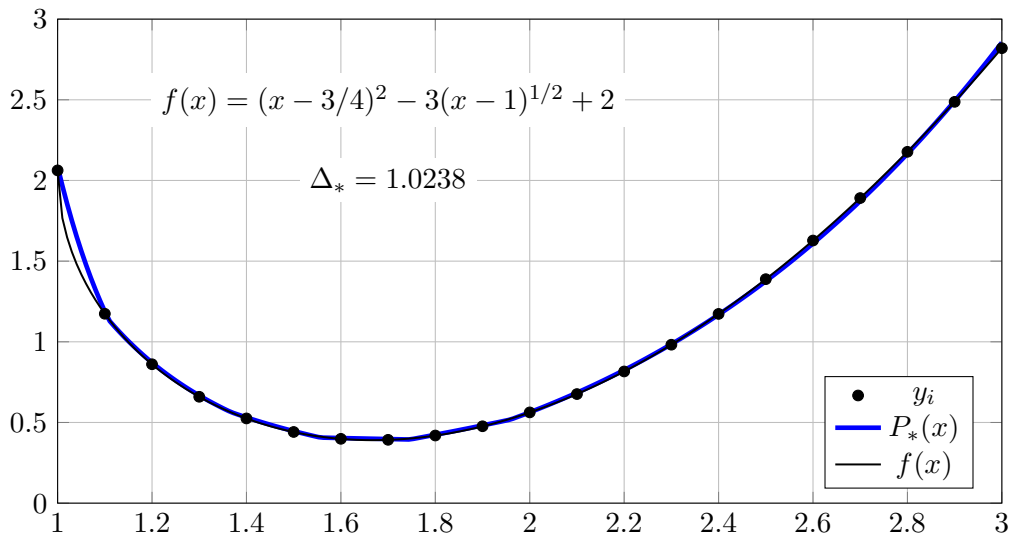

We conclude with results of approximation by max-algebra polynomials with $N=7$ monomials. The squared error obtained is $\Delta_{\ast}=1.0238$, whereas the approximating function (see Figure~\ref{F-AFMPN7}) is written in the ordinary notation as
\begin{multline*}
P_{\ast}(x)
=
\max(
2.0869x^{-5.9139},
1.5991x^{-3.3320},
1.2426x^{-2.5249},
\\
0.4555x^{-0.2581},
0.1055x^{2.3669},   
0.0351x^{4.0040},   
0.0351x^{4.0040}).
\end{multline*}

We observe that the last two monomials are identical, which means that the algorithm actually uses $N=6$ monomials to construct the optimal approximating polynomial rather than $7$. Moreover, as numerical experiments show, a further increase of the number $N$ of monomials does not improve the solution leaving the approximate polynomial unchanged.

\section{Conclusion}
\label{S-C}

In this paper, we have developed a new approach to solve discrete best approximation problems that are set in the framework of tropical algebra focused on the theory and applications of semirings and semifields with idempotent addition. Given sample points of the input and output of an unknown function defined on a tropical semifield, the problem is to approximate the function by generalized Puiseux polynomials in the sense of a generalized metric. The approximation problem involves evaluation of both unknown coefficients and exponents of monomials in the approximating polynomial. 

We first developed a solution approach to the problem defined in the max-plus algebra setting, and then extended the result to solve the problem in terms of max-algebra. The solution starts with the reduction of the polynomial approximation problem to the best approximation of a vector equation. This vector approximation problem has an analytical solution, which allows to represent both the unknown coefficients of monomials and approximation error in direct forms parameterized by the unknown exponents. As the next step of the approach, the problem of minimizing the approximation error is solved to obtain optimal values of the exponents. The obtained exponents are used to evaluate the coefficients, which together with the fixed exponents completely determine the best approximate polynomial.   

We rearranged the problem of minimization of the approximation error with respect to unknown exponents into combinatorial minimization over partitions of a finite set, where the evaluation of the objective function for each partition involves solving minimization problems that determine optimal exponents. To handle the combinatorial problem, we proposed a solution procedure based on the agglomerative clustering technique, which constructs an optimal (near-optimal) partition by successively merging clusters until a partition with a required number of parts is obtained. To solve the minimization problems for unknown exponents in the max-plus algebra setting, a transformation of the problems was implemented. This transformation turns the objective function into a polynomial where the unknown exponents become indeterminates, and hence allows to solve the problems as polynomial optimization problems that have known analytical solution. 

The solution developed in the framework of max-plus algebra can be directly applied to the solution of conventional problems of the Chebyshev discrete best approximation by piecewise linear functions. The solution in terms of max-algebra offers a way to solve approximation problems where approximating functions take the form of piecewise Puiseux polynomials. We evaluate the contribution of the study to be twofold. The proposed solutions provide an alternative approach to handle known approximation problems, and hence have strong potential as a new general tool to complement and supplement existing methods. At the same time, the results obtained further elaborate methods and extend the range of applications of tropical algebra, including solution of problems that involve uncertainty and incomplete data.

As a possible shortcoming of the obtained solutions, one can consider their time complexity that may be somewhat higher than that of other existing approximation methods. However, the computational complexity involved is still polynomial in both the number of sample points and the number of monomials in the approximating polynomial, and can be considered reasonable for many applications. Since the max-plus algebra based technique uses the approximation by convex piecewise linear functions, it cannot be applied to approximate arbitrary nonconvex functions, which presents another limitation of the approach. 

Possible areas of future investigation include the extension of the approximation approach to problems that are formulated in terms of other idempotent semifields than max-plus and max-algebra. Further development of the technique of finding exponents to ensure globally optimal solutions presents another promising research direction. Additional analysis of the computational complexity of the solution to improve the accuracy of the obtained estimates of complexity is also of interest.

\bibliographystyle{abbrvurl}

\bibliography{On_solution_of_tropical_discrete_best_approximation_problems}

\begin{thebibliography}{10}

\bibitem{Akian2011Best}
M.~Akian, S.~Gaubert, V.~Ni\c{t}ic\u{a}, and I.~Singer.
\newblock Best approximation in max-plus semimodules.
\newblock {\em Linear Algebra Appl.}, 435(12):3261--3296, 2011.
\newblock \href {https://doi.org/10.1016/j.laa.2011.06.009}
  {\path{doi:10.1016/j.laa.2011.06.009}}.

\bibitem{Anderberg1973Cluster}
M.~R. Anderberg.
\newblock {\em Cluster Analysis for Applications}.
\newblock Probability and Mathematical Statistics. Academic Press, New York,
  NY, 1973.
\newblock \href {https://doi.org/10.1016/C2013-0-06161-0}
  {\path{doi:10.1016/C2013-0-06161-0}}.

\bibitem{Butkovic2010Maxlinear}
P.~Butkovi\v{c}.
\newblock {\em Max-linear Systems}.
\newblock Springer Monographs in Mathematics. Springer, London, 2010.
\newblock \href {https://doi.org/10.1007/978-1-84996-299-5}
  {\path{doi:10.1007/978-1-84996-299-5}}.

\bibitem{Cameron1966Piecewise}
S.~H. Cameron.
\newblock Piece-wise linear approximation.
\newblock Technical Note CSTN-106, Computer Sciences Division, IIT Research
  Institute, Chicago, IL, 1966.

\bibitem{Camponogara2015Models}
E.~Camponogara and L.~F. Nazari.
\newblock Models and algorithms for optimal piecewise-linear function
  approximation.
\newblock {\em Math. Probl. Eng.}, 2015:876862, 2015.
\newblock \href {https://doi.org/10.1155/2015/876862}
  {\path{doi:10.1155/2015/876862}}.

\bibitem{Celikyilmaz2009Modeling}
A.~Celikyilmaz and I.~B. T\"{u}rksen.
\newblock {\em Modeling Uncertainty with Fuzzy Logic}, volume 240 of {\em
  Studies in Fuzziness and Soft Computing}.
\newblock Springer, Berlin, 2009.
\newblock \href {https://doi.org/10.1007/978-3-540-89924-2}
  {\path{doi:10.1007/978-3-540-89924-2}}.

\bibitem{Conn1988Computational}
A.~R. Conn and Y.~Li.
\newblock The computational structure and characterization of nonlinear
  discrete {C}hebyshev problem.
\newblock Technical Report 88-956, Department of Computer Science, Cornell
  University, Ithaca, NY, 1988.

\bibitem{Laplace1832Mecanique}
P.~S. de~Laplace.
\newblock {\em M{\'e}canique C{\'e}leste. Volume 2}.
\newblock Hillard, Gray, Littl{\`e}, and Wilkins, Boston, 1832.
\newblock (Engl. transl. with comment. by N.~Bowditch).

\bibitem{Dubois1980Fuzzy}
D.~J. Dubois and H.~M. Prade.
\newblock {\em Fuzzy Sets and Systems}, volume 144 of {\em Mathematics in
  Science and Engineering}.
\newblock Academic Press, San Diego, 1980.

\bibitem{Esparza2008Approximative}
J.~Esparza, T.~Gawlitza, S.~Kiefer, and H.~Seidl.
\newblock Approximative methods for monotone systems of min-max-polynomial
  equations.
\newblock In L.~Aceto, I.~Damg{\aa}rd, L.~A. Goldberg, M.~M. Halld{\'o}rsson,
  A.~Ing{\'o}lfsd{\'o}ttir, and I.~Walukiewicz, editors, {\em Automata,
  Languages and Programming}, volume 5125 of {\em Lecture Notes in Computer
  Science}, pages 698--710, Berlin, 2008. Springer.
\newblock \href {https://doi.org/10.1007/978-3-540-70575-8\_57}
  {\path{doi:10.1007/978-3-540-70575-8\_57}}.

\bibitem{Gluss1962Further}
B.~Gluss.
\newblock Further remarks on line segment curve-fitting using dynamic
  programming.
\newblock {\em Commun. ACM}, 5(8):441--443, 1962.
\newblock \href {https://doi.org/10.1145/368637.368753}
  {\path{doi:10.1145/368637.368753}}.

\bibitem{Golan2003Semirings}
J.~S. Golan.
\newblock {\em Semirings and Affine Equations Over Them}, volume 556 of {\em
  Mathematics and Its Applications}.
\newblock Springer, Dordrecht, 2003.
\newblock \href {https://doi.org/10.1007/978-94-017-0383-3}
  {\path{doi:10.1007/978-94-017-0383-3}}.

\bibitem{Gondran2008Graphs}
M.~Gondran and M.~Minoux.
\newblock {\em Graphs, Dioids and Semirings}, volume~41 of {\em Operations
  Research/ Computer Science Interfaces}.
\newblock Springer, New York, NY, 2008.
\newblock \href {https://doi.org/10.1007/978-0-387-75450-5}
  {\path{doi:10.1007/978-0-387-75450-5}}.

\bibitem{Grigoriev2018Tropical}
D.~Grigoriev.
\newblock Tropical {N}ewton-{P}uiseux polynomials.
\newblock In V.~P. Gerdt, W.~Koepf, W.~M. Seiler, and E.~V. Vorozhtsov,
  editors, {\em Computer Algebra in Scientific Computing}, volume 11077 of {\em
  Lecture Notes in Computer Science}, pages 177--186. Springer, Cham, 2018.
\newblock \href {https://doi.org/10.1007/978-3-319-99639-4\_12}
  {\path{doi:10.1007/978-3-319-99639-4\_12}}.

\bibitem{Grigoriev2014Tropical}
D.~Grigoriev and V.~Shpilrain.
\newblock Tropical cryptography.
\newblock {\em Comm. Algebra}, 42(6):2624--2632, 2014.
\newblock \href {https://doi.org/10.1080/00927872.2013.766827}
  {\path{doi:10.1080/00927872.2013.766827}}.

\bibitem{Grigoriev2019Tropical}
D.~Grigoriev and V.~Shpilrain.
\newblock Tropical cryptography {II}: Extensions by homomorphisms.
\newblock {\em Comm. Algebra}, 47(10):4224--4229, 2019.
\newblock \href {https://doi.org/10.1080/00927872.2019.1581213}
  {\path{doi:10.1080/00927872.2019.1581213}}.

\bibitem{Heidergott2006Maxplus}
B.~Heidergott, G.~J. Olsder, and J.~{van der Woude}.
\newblock {\em Max {P}lus at Work}.
\newblock Princeton Series in Applied Mathematics. Princeton Univ. Press,
  Princeton, NJ, 2006.

\bibitem{Imai1986Optimal}
H.~Imai and M.~Iri.
\newblock An optimal algorithm for approximating a piecewise linear function.
\newblock {\em J. Inf. Process.}, 9(3):159--162, 1986.

\bibitem{Itenberg2007Tropical}
I.~Itenberg, G.~Mikhalkin, and E.~Shustin.
\newblock {\em Tropical Algebraic Geometry}, volume~35 of {\em Oberwolfach
  Seminars}.
\newblock Birkh\"{a}user, Basel, 2007.
\newblock \href {https://doi.org/10.1007/978-3-7643-8310-7}
  {\path{doi:10.1007/978-3-7643-8310-7}}.

\bibitem{Krivulin2012Solution}
N.~Krivulin.
\newblock A solution of a tropical linear vector equation.
\newblock In S.~Yenuri, editor, {\em Advances in Computer Science}, volume~5 of
  {\em Recent Advances in Computer Engineering Series}, pages 244--249. WSEAS
  Press, Athens, 2012.
\newblock \href {https://arxiv.org/abs/1212.6107} {\path{arXiv:1212.6107}}.

\bibitem{Krivulin2013Solution-linear}
N.~Krivulin.
\newblock Solution of linear equations and inequalities in idempotent vector
  spaces.
\newblock {\em Int. J. Appl. Math. Inform.}, 7(1):14--23, 2013.
\newblock \href {https://arxiv.org/abs/1305.4300} {\path{arXiv:1305.4300}}.

\bibitem{Krivulin2021Algebraic}
N.~Krivulin.
\newblock Algebraic solution of tropical polynomial optimization problems.
\newblock {\em Mathematics}, 9(19):2472, 2021.
\newblock \href {https://arxiv.org/abs/2011.05460} {\path{arXiv:2011.05460}},
  \href {https://doi.org/10.3390/math9192472} {\path{doi:10.3390/math9192472}}.

\bibitem{Krivulin2023Algebraic}
N.~Krivulin.
\newblock Algebraic solution of tropical best approximation problems.
\newblock {\em Mathematics}, 11(18):3949, 2023.
\newblock \href {https://arxiv.org/abs/2308.07210} {\path{arXiv:2308.07210}},
  \href {https://doi.org/10.3390/math11183949}
  {\path{doi:10.3390/math11183949}}.

\bibitem{Krivulin2009Solution}
N.~K. Krivulin.
\newblock On solution of a class of linear vector equations in idempotent
  algebra.
\newblock {\em Vestnik Sankt-Peterburgskogo Universiteta. Seriya~10},
  (3):64--77, 2009.
\newblock (in Russian).

\bibitem{Li1992Morphological}
D.~Li.
\newblock Morphological template decomposition with max-polynomials.
\newblock {\em J. Math. Imaging Vision}, 1(3):215--221, 1992.
\newblock \href {https://doi.org/10.1007/BF00129876}
  {\path{doi:10.1007/BF00129876}}.

\bibitem{Maclagan2015Introduction}
D.~Maclagan and B.~Sturmfels.
\newblock {\em Introduction to Tropical Geometry}, volume 161 of {\em Graduate
  Studies in Mathematics}.
\newblock AMS, Providence, RI, 2015.
\newblock \href {https://doi.org/10.1090/gsm/161} {\path{doi:10.1090/gsm/161}}.

\bibitem{Maragos2021Tropical}
P.~Maragos, V.~Charisopoulos, and E.~Theodosis.
\newblock Tropical geometry and machine learning.
\newblock {\em Proc. IEEE}, 109(5):728--755, 2021.
\newblock \href {https://doi.org/10.1109/JPROC.2021.3065238}
  {\path{doi:10.1109/JPROC.2021.3065238}}.

\bibitem{Markwig2010Field}
T.~Markwig.
\newblock A field of generalised {P}uiseux series for tropical geometry.
\newblock {\em Rend. Sem. Mat. Univ. Politec. Torino}, 68(1):79--82, 2010.

\bibitem{Mhaskar2000Fundamentals}
H.~N. Mhaskar and D.~V. Pai.
\newblock {\em Fundamentals of Approximation Theory}.
\newblock Narosa Publishing House, New Delhi, 2000.

\bibitem{Osborne1967Best}
M.~R. Osborne and G.~A. Watson.
\newblock On the best linear {C}hebyshev approximation.
\newblock {\em Comput. J.}, 10(2):172--177, 1967.
\newblock \href {https://doi.org/10.1093/comjnl/10.2.172}
  {\path{doi:10.1093/comjnl/10.2.172}}.

\bibitem{Saad2021Zerosum}
O.~Saadi.
\newblock {\em Zero-sum repeated games: Accelerated algorithms and tropical
  best-approximation}.
\newblock PhD thesis, Institut Polytechnique de Paris, 2021.

\bibitem{Samarasinghe2006Neural}
S.~Samarasinghe.
\newblock {\em Neural Networks for Applied Sciences and Engineering}.
\newblock Auerbach Publ., New York, 2006.
\newblock \href {https://doi.org/10.1201/9780849333750}
  {\path{doi:10.1201/9780849333750}}.

\bibitem{Schunn2010Uncertainly}
C.~D. Schunn.
\newblock From uncertainly exact to certainly vague: Epistemic uncertainty and
  approximation in science and engineering problem solving.
\newblock In B.~H. Ross, editor, {\em The Psychology of Learning and
  Motivation}, volume~53 of {\em Psychology of Learning and Motivation}, pages
  227--252. Academic Press, San Diego, CA, 2010.
\newblock \href {https://doi.org/10.1016/S0079-7421(10)53006-8}
  {\path{doi:10.1016/S0079-7421(10)53006-8}}.

\bibitem{Steffens2006History}
K.-G. Steffens.
\newblock {\em The History of Approximation Theory}.
\newblock Birkh\"{a}user, Boston, MA, 2006.
\newblock \href {https://doi.org/10.1007/0-8176-4475-X}
  {\path{doi:10.1007/0-8176-4475-X}}.

\bibitem{Stone1961Approximation}
H.~Stone.
\newblock Approximation of curves by line segments.
\newblock {\em Math. Comp.}, 15:40--47, 1961.
\newblock \href {https://doi.org/10.1090/S0025-5718-1961-0119390-6}
  {\path{doi:10.1090/S0025-5718-1961-0119390-6}}.

\bibitem{Szusz2010Linear}
E.~K. Szusz and A.~R. Willms.
\newblock A linear time algorithm for near minimax continuous piecewise linear
  representations of discrete data.
\newblock {\em SIAM J. Sci. Comput.}, 32(5):2584--2602, 2010.
\newblock \href {https://doi.org/10.1137/090769077}
  {\path{doi:10.1137/090769077}}.

\bibitem{Tharwat2010One}
A.~Tharwat and K.~Zimmermann.
\newblock One class of separable optimization problems: Solution method,
  application.
\newblock {\em Optimization}, 59(5):619--625, 2010.
\newblock \href {https://doi.org/10.1080/02331930801954698}
  {\path{doi:10.1080/02331930801954698}}.

\bibitem{Tomek1974Two}
I.~Tomek.
\newblock Two algorithms for piecewise-linear continuous approximation of
  functions of one variable.
\newblock {\em IEEE Trans. Comput.}, 23(4):445--448, 1974.
\newblock \href {https://doi.org/10.1109/T-C.1974.223961}
  {\path{doi:10.1109/T-C.1974.223961}}.

\bibitem{Ward1963Hierarchical}
J.~H. {Ward Jr.}
\newblock Hierarchical grouping to optimize an objective function.
\newblock {\em J. Amer. Statist. Assoc.}, 58(301):236--244, 1963.
\newblock \href {https://doi.org/10.1080/01621459.1963.10500845}
  {\path{doi:10.1080/01621459.1963.10500845}}.

\bibitem{Watson1970Algorithm}
G.~A. Watson.
\newblock On an algorithm for nonlinear minimax approximation.
\newblock {\em Commun. ACM}, 13(3):160--162, 1970.
\newblock \href {https://doi.org/10.1145/362052.362056}
  {\path{doi:10.1145/362052.362056}}.

\bibitem{Zhang2018Tropical}
L.~Zhang, G.~Naitzat, and L.-H. Lim.
\newblock Tropical geometry of deep neural networks.
\newblock In J.~Dy and A.~Krause, editors, {\em Proceedings of the 35th
  International Conference on Machine Learning}, volume~80 of {\em Proceedings
  of Machine Learning Research}, pages 5824--5832. PMLR, Cambridge, MA, 2018.

\bibitem{Zimmermann1984Maxseparable}
K.~Zimmermann.
\newblock On max-separable optimization problems.
\newblock In R.~E. Burkard, R.~A. Cuninghame-Green, and U.~Zimmermann, editors,
  {\em Algebraic and Combinatorial Methods in Operations Research}, volume~95
  of {\em North-Holland Mathematics Studies}, pages 357--362. North-Holland,
  Amsterdam, 1984.
\newblock \href {https://doi.org/10.1016/S0304-0208(08)72967-0}
  {\path{doi:10.1016/S0304-0208(08)72967-0}}.

\bibitem{Zimmermann2003Disjunctive}
K.~Zimmermann.
\newblock Disjunctive optimization, max-separable problems and extremal
  algebras.
\newblock {\em Theoret. Comput. Sci.}, 293(1):45--54, 2003.
\newblock \href {https://doi.org/10.1016/S0304-3975(02)00231-1}
  {\path{doi:10.1016/S0304-3975(02)00231-1}}.

\end{thebibliography}

\end{document}